\documentclass{amsart}
\usepackage{amsfonts,amssymb,stmaryrd,amscd,amsmath,latexsym,amsbsy}
\numberwithin{equation}{section}
\theoremstyle{plain}
\newtheorem{thm}{Theorem}[section]

\newtheorem{prop}[thm]{Proposition}

\newtheorem{defi}[thm]{Definition}
\newtheorem{lem}[thm]{Lemma}
\newtheorem{cor}[thm]{Corollary}

\theoremstyle{remark}
\newtheorem{rema}[thm]{Remark}

\renewcommand{\_}[1]{_{(#1)}}

\newcommand{\ad}{{\mbox{\upshape{ad}}}}

\newcommand{\al }{\alpha }

\newcommand{\C}{{\mathbb C}}
\newcommand{\Char}{\mathrm{Char}}

\newcommand{\N}{{\mathbb N}}

\newcommand{\cF}{{\mathcal F}}

\newcommand{\cqg}{{\qfield_q[G]}}

\newcommand{\field}{k}

\newcommand{\gfrak}{{\mathfrak g}}

\newcommand{\Hom}{{\mathrm{Hom}}}

\newcommand{\hght}{\mathrm{ht}}

\newcommand{\id}{{\mbox{id}}}

\newcommand{\kow}{{\varDelta}}

\newcommand{\Lin}{\mbox{span}}

\newcommand{\Map}{\mathrm{Map}}

\newcommand{\ot}{\otimes}

\newcommand{\pr}{\mathrm{pr}}

\newcommand{\qfield}{\field }

\newcommand{\rcs }{{\mathfrak K}}

\newcommand{\rk}{\mathrm{rank}}

\newcommand{\rank}{\mathrm{rank}}

\newcommand{\slfrak}{{\mathfrak{sl}}}
\newcommand{\spec}{\mathrm{spec}}
\newcommand{\supp}{\mathrm{supp}}

\newcommand{\Uq}{U}

\newcommand{\uqbp}{U^{\ge 0}} 

\newcommand{\uqg}{{U_q(\mathfrak{g})}}

\newcommand{\vep}{\varepsilon}

\newcommand{\weights}{\Lambda}

\newcommand{\wurz}{\Phi}

\newcommand{\Z}{{\mathbb Z}}

\begin{document}
\title[Right coideal subalgebras of $U^{\ge 0}$]{Right coideal subalgebras
of the Borel part of a quantized enveloping algebra}

\author{ Istv{\'a}n Heckenberger}
\address{ Istv{\'a}n Heckenberger, Fachbereich Mathematik, Philipps-Universit{\"a}t Marburg, Hans-Meerwein-Stra{\ss}e, 35032 Marburg, Germany}
\email{heckenberger@mathematik.uni-marburg.de}

\author{Stefan Kolb}
\address{Stefan Kolb, Korteweg-de Vries Institute for Mathematics,
University of Amsterdam, Science Park 904, 1098 XH Amsterdam,
The Netherlands}

\email{s.kolb@uva.nl}

\thanks{The work of I.~Heckenberger was supported by DFG within a Heisenberg
fellowship. The work of S.~Kolb was supported by the School of Mathematics and
the Maxwell Institute for Mathematical Sciences at the University of Edinburgh
and by the Netherlands Organization for Scientific Research (NWO) within the
VIDI-project ``Symmetry and modularity in exactly solvable models''.}


\begin{abstract}
  For the Borel part of a quantized enveloping algebra we classify all
  right coideal subalgebras for which the intersection with the
  coradical is a Hopf algebra. The result is expressed in terms of
  characters of the subalgebras $U^+[w]$ of the
  quantized enveloping algebra, introduced by de Concini, Kac, and
  Procesi for any Weyl group element $w$.
  We explicitly determine all characters
  of $U^+[w]$ building on recent work by Yakimov on prime ideals of
  $U^+[w]$ which are invariant under a torus action.
\end{abstract}

\maketitle

\section{Introduction}
Let $\gfrak$ be a complex, finite-dimensional, semisimple Lie algebra
and $\uqg$ the corresponding quantized enveloping algebra. Quantum
analogs of Lie subalgebras of $\gfrak$ are often 
realized as coideal subalgebras of $\uqg$. Recall that a subalgebra $C$
of a Hopf algebra $H$ is called a right coideal subalgebra if the
coproduct $\kow$ of $H$ satisfies $\kow(C)\subseteq C\ot H$. The
universal enveloping algebra $U(\gfrak)$ is a cocommutative Hopf
algebra and hence right coideal subalgebras of $U(\gfrak)$ are always Hopf
subalgebras. The quantum deformation $\uqg$, however, is essentially
obtained by deforming the coproduct. Hence one expects quantum
deformations of Hopf subalgebras of $U(\gfrak)$ inside $\uqg$ only to
satisfy the weaker coideal property. 

Let $\Pi$ be a basis of simple roots for $\gfrak$. Let $U=\uqg$ be
defined as in \cite[Chapter 4]{b-Jantzen96} over a base field
$\field$ with $q\in \field\setminus \{0\}$ not a root
of unity. Let $U^+$ and $U^0$ denote the subalgebras of $U$ generated
by the sets $\{E_\alpha\,|\,\alpha\in \Pi\}$ and $\{K_\alpha, 
K_\alpha^{-1}\,|\,\alpha\in \Pi\}$, respectively, and define the
positive Borel part $U^{\ge 0}= U^+U^0$ which is a Hopf subalgebra of
$U$. In the present paper we give an explicit combinatorial
classification of all right coideal subalgebras $C$ of $U^{\ge 0}$ for which
$C\cap U^0$ is a Hopf algebra. To be more precise, let $Q$ denote the
root lattice and, for any element $w$ in the Weyl group $W$ of $\gfrak$,
let $U^+[w]$ denote the subalgebra of $U^+$ generated by the
corresponding Lusztig root vectors as defined in
\cite[8.24]{b-Jantzen96}. Our first main result, which is made precise
in Theorem \ref{thm:coid-class}, states the following:

\smallskip

\noindent \textit{There is a canonical bijection between the set of
  all right coideal subalgebras $C$ of $U^{\ge 0}$ for which $C\cap U^0$
  is a Hopf algebra and the set of all triples $(w,\phi,L)$ where $w\in
  W$, $\phi:U^+[w]\rightarrow \field$ is a character, and $L$ is a
  subgroup of $Q$ such that $\phi$ and $L$ satisfy an additional compatibility
  condition.}

\smallskip

Many examples of coideal subalgebras of $U$ were initially constructed
with a theory of homogeneous spaces or harmonic analysis for quantum
groups in mind.  Koornwinder's observation \cite{a-Koornwinder93} that
a quantum group analog of one-dimensional complex projective space
can be obtained via a skew primitive element in
$U_q(\slfrak_2(\C))$ inspired the development of a theory of quantum
symmetric spaces via coideal subalgebras. First, quantum analogs of
all classical symmetric pairs were constructed in a case by case
fashion \cite{a-Noumi96}, \cite{a-Dijk96}, \cite{a-NDS97}. Later, a
comprehensive theory of quantum symmetric pairs as one-sided coideal subalgebras
in $U$ was developed by G.~Letzter \cite{MSRI-Letzter},
\cite{a-Letzter03}. Letzter's work also contains general qualitative
results about the structure of coideal subalgebras of $U$, mainly in
terms of filtrations and the associated graded algebras \cite[Section
4]{MSRI-Letzter}. Closely related is a program initiated by Kharchenko
to determine all right coideal subalgebras of $U^{\ge 0}$ which
contain $U^0$. It was proved for $\gfrak$ of type $A_n$
\cite{a-KharSag08}, $B_n$ \cite{a-Khar09p}, and $G_2$
\cite{a-Pogorelsky09} that the number of such right coideal subalgebras
coincides with the order of the Weyl group. Recently, the situation was
clarified by Schneider and the first named author who proved in
\cite{a-HeckSchn09p} that the algebras $U^+[w]U^0$, where $w\in W$, exhaust
all right coideal subalgebras of $U^{\ge 0}$ which contain $U^0$.
This result forms the starting point of
the present work. 

Motivated by the above classification, the second main result of this
paper consists of an explicit combinatorial description of the set
$\Char(U^+[w])$ of all characters of $U^+[w]$. Define
$W^{\le w}= \{y\in W\,|\,y\le w\}$ where $\le$ denotes the
Bruhat order on $W$.
In Sections \ref{sec:H-strat} and \ref{sec:char-class} we show
the following:   

\smallskip
\noindent
\textit{The set $\Char(U^+[w])$ can be canonically identified with a disjoint
union of spectra of Laurent polynomial rings indexed by the elements
of a subset $W^w\subseteq W^{\le w}$. The subset  
$W^w$ is given explicitly, and the dimension of the component
corresponding to $y\in W^w$ is $\ell(w)-\ell(y)$, where $\ell$ denotes
the length function.}
\smallskip

To determine $\Char(U^+[w])$ we apply the theory of $H$-stable prime
ideals, outlined in \cite[Part II]{b-BG02}. Here
$H=(\field\setminus\{0\})^{\rank(\gfrak)}$ is a torus which naturally
acts on $U^+[w]$. A subspace of $U^+[w]$ is $H$-stable if and only
if it is naturally graded by the root lattice $Q$. Following a standard
construction, one associates an $H$-prime ideal of $U^+[w]$ to any
character of $U^+[w]$. Recently, Yakimov obtained an order 
preserving bijection between $W^{\le w}$ with the Bruhat order and the
set of $H$-prime ideals of $U^+[w]$, which is ordered by inclusion
\cite{a-Yakimov09p}. To determine $\Char(U^+[w])$ we directly find a large
family of $H$-prime ideals which correspond to characters.
We then apply Yakimov's result and previous results by Gorelik
\cite{a-Gorelik00} to show that this family is complete. 

Characters of $U^+[w]$ appear in the classification of right coideal
subalgebras via a standard construction: For any right coideal subalgebra
$C$ of a Hopf algebra $H$ and any character $\phi:C\rightarrow \field$
one obtains a new right coideal subalgebra $C_\phi$ by application of the
coproduct and evaluation of $\phi$ on the first tensor component
(cf.~Section \ref{sec:red-char}). It was recently noted that this
construction is also at the heart of coideal subalgebras of $U$ which
appear in the theory of quantum symmetric pairs \cite[Section
4]{a-KS09p}. One  might ask if $U$ possesses certain standard right
coideal subalgebras from which all other right coideal subalgebras are
obtained via characters in this way. This is but one question on the
way towards a general combinatorial classification of right coideal
subalgebras of $U$.

The classification of right coideal subalgebras of $\uqbp $ which contain
$U^0$ also holds for small quantum groups of semisimple Lie algebras
where $q$ is a root of unity \cite{a-KharSag08},
\cite{a-HeckSchn09p}. In the setting of the present paper, however,
the assumption that $q$ is not a root of unity is essential. It would 
be interesting to extend the classification of right coideal
subalgebras to more general classes of pointed Hopf algebras.  

This paper consists of the two sections outlined above and an
appendix. The appendix contains technical results on root systems and
Weyl group combinatorics which are used to prove the statements
leading up to Theorem \ref{thm:char-class}.   

\section{Right coideal subalgebras of $\uqbp$}
\subsection{Quantized enveloping algebras and right coideal subalgebras}
\label{sec:qea}

We mostly follow the notation and conventions of \cite{b-Jantzen96}.
Let $\gfrak $ be a finite-dimensional complex semisimple Lie algebra and let
$\wurz $ be the root system with respect to a fixed Cartan
subalgebra. We also fix a basis $\Pi$ of  $\Phi$ and denote by
$\wurz^+$ and $\wurz^-$ the corresponding sets of positive roots and
negative roots, respectively.
Let $W$ be the Weyl group of $\gfrak $ and let $(\cdot ,\cdot )$ be
the invariant scalar product on the real vector space generated by
$\Pi $ such that $(\al ,\al )=2$ for all short roots in each
component. For any $\beta\in \wurz$ we write $s_\beta$ to denote the
reflection at the hyperplane orthogonal to $\beta$ with respect to
$(\cdot, \cdot)$. Let $Q=\Z \Pi $ be the root lattice and let $Q_+=\N _0\Pi
$. For each $\al \in \Pi $ let $d_\al =(\al ,\al )/2$.
Let $U=U_q(\gfrak )$ be the quantized enveloping algebra of $\gfrak $
in the sense of \cite[Chapter 4]{b-Jantzen96}.
More precisely, let $\field $ be a field and
fix an element $q\in \field $ with $q\not=0$ and $q^n\not=1$ for all
$n\in \N $. Then $U$ is the
unital associative algebra defined over $\field $
with generators $K_\al , K_\al ^{-1}, E_\al , F_\al $ for all $\al \in \Pi
$ and relations given in \cite[4.3]{b-Jantzen96}. 
By \cite[Proposition 4.11]{b-Jantzen96}
there is a unique Hopf algebra structure on $U$ with coproduct $\kow$,
counit $\vep$, and antipode $S$ such that
\begin{align}
  \label{eq:UHopf1}
  \kow (E_\al )=&E_\al \ot 1+K_\al \ot E_\al ,&
  \vep (E_\al )=&0,& S(E_\al)&=-K_\al^{-1} E_\al ,\\
  \kow (F_\al )=&F_\al \ot K_\al ^{-1}+1 \ot F_\al ,&
  \vep (F_\al )=&0,& S(F_\al)&=-F_\al K_\al,\\
  \kow (K_\al )=&K_\al \ot K_\al ,& \vep (K_\al )=&1, &S(K_\al)&=K_\al^{-1}.
  \label{eq:UHopf3}
\end{align}
We will make use of Sweedler notation for the coproduct in the form
$\kow(x)=x_{(1)}\ot x_{(2)}$ for any $x\in U$.
Let $\ad $ denote the left adjoint action of $U$ on itself, that is,
$(\ad x)(y)=x\_1 y S(x\_2)$ for all $x,y\in U$.

As in \cite[Chapter 4]{b-Jantzen96} let $U^+$, $U^0$, and $\uqbp $
be the subalgebras of $U$ generated by the sets $\{E_\al \,|\,\al \in \Pi
\}$, $\{K_\al ,K_\al ^{-1}\,|\,\al \in \Pi \}$
and $\{ K_\al ,K_\al ^{-1},E_\al \,|\,\al \in \Pi \}$, respectively.
In fact, $U^0$ and $\uqbp $ are Hopf subalgebras of $\uqg $.
For any $\beta \in Q$ define $K_\beta =\prod _{\al \in
\Pi }K_\al ^{n_\al }$ if $\beta =\sum _{\al \in \Pi }n_\al \al $ for
some $n_\alpha\in \Z$.

For any $\al \in \Pi$ let $T_\al $ denote the algebra automorphisms of $U$
defined in \cite[8.14]{b-Jantzen96}.
Let $w\in W$ be an element of length $\ell (w)=t$ and choose $\al
_1,\dots,\al _t\in \Pi $ such that $s_{\al _1}s_{\al _2}\cdots s_{\al
  _t}$ is a reduced expression of $w$. For all $i\in \{1,2,\dots,t\}$
let $\beta _i=s_{\al _1}\cdots s_{\al _{i-1}}\al _i$.
In \cite[8.21]{b-Jantzen96} Lusztig root
vectors in $U^+$ are defined by $E_{\beta _i}=T_{\al_1}\cdots
T_{\al_{i-1}}E_{\al_i}$ for all  $i\in \{1,\dots,t\}$. Following \cite{a-DCoKacPro93} the subspace 
\[ U^+[w]=\Lin _\field \{ E_{\beta _t}^{a_t} \cdots E_{\beta _2}^{a_2}
  E_{\beta _1}^{a_1} \,|\,a_1,\dots,a_t\in \N _0 \}, \]
is attached to $w$
in \cite[8.24]{b-Jantzen96}.
It is shown in \cite{a-DCoKacPro93} that $U^+[w]$ is a subalgebra
which does not depend on the reduced expression. The following
observation is the starting point of our investigations of right
coideal subalgebras of $\uqbp $. 

\begin{thm} \cite[Theorem 7.3]{a-HeckSchn09p} \label{th:HS}
  The map from $W$ to the set of right coideal subalgebras of
  $\uqbp $ containing $U^0$, given by $w\mapsto U^+[w]U^0$, is
  a bijection.
\end{thm}

For all $\al \in Q_+$ let
\begin{align}
  U^+_\al =&\{x\in U^+\,|\,K_\beta x K_\beta ^{-1}=
  q^{(\beta ,\al )}x \quad \text{for all $\beta \in Q$}\}.
\end{align}
The algebra $\uqbp $ admits a $Q^2$-grading given by
\begin{align}
  \label{eq:Z2grading}
  \uqbp =&\mathop{\oplus }_{\alpha ,\beta \in Q}U^+_\alpha K_{\beta }.
\end{align}
For all $\alpha ,\beta \in Q$
let $\pr _{(\al ,\beta )}:\uqbp \to U^+_\al K_\beta $
be the unique $Q^2$-graded projection. Note that $\pr _{(\al ,\beta )}=0$
if $\al \in Q\setminus Q_+$.

A subspace $T$ of $U^0$ is a Hopf subalgebra if and only if there
exists a subgroup $L$ of the abelian group $Q$ such that $T=\Lin
_\field \{K_\al \,|\,\alpha \in L\}$. We write $T_L$ for this Hopf
algebra. In particular, $T_{\{0\}}=\field$ and $T_Q=U^0$. 
We will frequently use the notation
\begin{align*}
  L^\perp =\{\gamma \in Q\,|\,(\gamma ,\beta)=0 \mbox{ for all
    $\beta\in L$}\}
\end{align*}
to denote the orthogonal complement of a subgroup $L\subseteq
Q$. The adjoint action of $T_L$ on $\uqbp$ is diagonalizable.
We say that an element $x\in \uqbp \setminus \{0\}$ is
a \textit{weight vector for}
$\ad \,T_L$,
if there exists $\al \in Q$ such that
\begin{align}
  (\ad \,K_\beta )(x)=q^{(\beta ,\al )}x \quad \text{for all $\beta \in L$.}
  \label{eq:adTLw}
\end{align}
Any $\ad\,T_L$-stable subspace of $\uqbp $ has a basis consisting of
weight vectors for $\ad \,T_L$. 

\begin{lem} \label{le:adTLw}
  Let $L\subseteq Q$ be a subgroup and $\al _1,\al _2,\beta
  _1,\beta _2\in Q$. Let $x\in \uqbp $ be a weight vector
  for $\ad \,T_L$ such that $\pr _{(\al _1,\beta _1)}(x)\not=0$ and
  $\pr _{(\al _2,\beta _2)}(x)\not=0$. Then $\al _1-\al _2\in L^\perp $.
\end{lem}

\begin{proof}
  The assumptions of the Lemma together with the direct sum
  decomposition \eqref{eq:Z2grading} imply that $(\ad
  K_\beta)(x)=q^{(\alpha_1,\beta)}x=q^{(\alpha_2,\beta)}x $ for all
  $\beta\in L$. As $q$ is not a root of unity one obtains
  $(\alpha_1-\alpha_2,\beta)=0$ for all $\beta\in L$.
\end{proof}

For any right coideal $C\subseteq \uqbp$ and any $\beta\in Q$ define
\begin{align}\label{eq:Cbeta}
  C_\beta=U^+K_\beta\cap C.
\end{align}
The following lemma is an adapted version of
\cite[Lemmata\,1.1,\,1.3]{MSRI-Letzter}.

\begin{lem}\label{le:Qdeg}
  Let $C\subseteq \uqbp $ be a right coideal. Then $C=\oplus _{\beta
    \in Q}C_\beta$. If $C$ is an algebra then this decomposition is an
  algebra grading of $C$ by $Q$.
\end{lem}

\begin{proof}
  Let $p:\uqbp \to U^0$, $p=\oplus _{\beta \in Q}\pr _{(0,\beta )}$,
  be the unique $Q^2$-graded projection. As $p$ is an algebra map,
  Equations~\eqref{eq:UHopf1}
  and \eqref{eq:UHopf3} imply that $p$ is a coalgebra homomorphism.
  Since $C$ is a right coideal,
  the map $p'=(\id \ot p)\kow :\uqbp \to \uqbp \ot U^0$
  induces a right coaction of $U^0$ on $C$. Further,
  \begin{align*}
    p'(x K_\beta)=x K_\beta \ot K_\beta \quad\mbox{for all }
    \beta\in Q, \quad x\in U^+,
  \end{align*}
  and hence $C=\oplus _{\beta \in Q}C_\beta $.
  The last claim of the lemma follows now from the commutation
  relations in $\uqbp$.
\end{proof}
The comultiplication $\kow $ of $\uqbp $
is compatible with the $Q^2$-grading \eqref{eq:Z2grading} via
\begin{align}
  \kow (x)-x\ot K_\beta 
  \in \mathop{\oplus }_{\gamma \prec \al }U^+_\gamma
  K_{\al +\beta -\gamma }
  \otimes U^+_{\al -\gamma } K_\beta 
  \label{eq:Q2kow}
\end{align}
for $x\in U^+_\al K_\beta $, $\alpha ,\beta \in Q$, where $\gamma \prec \al $
means that $\al -\gamma \in Q_+\setminus \{0\}$.
The relation $\prec $ defines a partial ordering on $Q$.

The next lemmata show that left coideals are useful for the study of
right coideals. 

\begin{lem}\label{le:tgamma}
  Let $C\subseteq \uqbp $ be a right coideal. Let $\beta \in Q$ and
  $x\in C_\beta \setminus \{0\}$.
  Let $J$ be a $Q^2$-homogeneous subspace and a left coideal of $\uqbp $.
  Let $\gamma \in Q$ be maximal with respect to $\prec $ such that
  $\pr _{(\gamma ,\beta )}(x)\notin J$. Then $K_{\beta +\gamma }\in C$.
\end{lem}

\begin{proof}
  For all $\al \in Q_+$ let $x_\al \in U^+_\al $ be
  such that $x=\sum _{\al \in Q_+}x_\alpha K_\beta $.
  Let $p_J:\uqbp \to \uqbp /J$ be the canonical map of left $\uqbp $-comodules.
  Since $J$ is $Q^2$-homogeneous, the map $p_J$ is $Q^2$-graded.
  Let $\alpha \in Q_+$ with $\gamma \prec \alpha $.
  By assumption, $x_\al K_\beta \in J$.
  Since $J$ is a left coideal, we conclude that
  $(\id \ot p_J)\kow (x_\alpha K_\beta )=0$.
  Thus the maximality of $\gamma $ and Relation~\eqref{eq:Q2kow} imply that
  \begin{align*}
    (\id \ot p_J)\kow (x)=&K_{\beta +\gamma }\ot p_J(x_\gamma K_\beta)\\
    &+ \text{terms $x_1\otimes p_J(x_2)$ with $x_2\in U^+_\al
      K_\beta$ for some $\al \neq \gamma$.}
  \end{align*}
  Since $C$ is a right coideal one obtains $K_{\beta +\gamma }\in C$.
\end{proof}

\begin{lem}
  \label{le:wvec}
  Let $C\subseteq \uqbp $ be a right coideal and let $x\in C_{-\beta
  }$ for some $\beta \in Q$. 
  Let $L\subseteq Q$ be a subgroup such that $C\cap U^0\subseteq T_L$.
  Assume that $x$ is a weight vector for $\ad \,T_L$. Let $\al\in Q_+$
  be maximal with respect to $\prec$ such that $x_\al:=\pr
  _{(\al,-\beta )}(x)\not=0$. Let $J_\al \subseteq \uqbp $
  be the left coideal generated by $x_\alpha$.

  (1) If $\gamma \in Q$ satisfies $x_\gamma :=\pr _{(\gamma
    ,-\beta)}(x) \neq 0$ then $x_\gamma \in J_\al$ and $\al -\gamma \in
  L^\perp \cap Q_+$.

  (2) The element $\alpha\in Q_+$ is uniquely determined by $x$.
  
  (3) If $\pr _{(\gamma ,-\beta )}(x)\not=0$ for some $\gamma \in \beta
  +L^\perp $ then $\beta =\al $.
\end{lem}

\begin{proof}
  Since $x_\al$ is $Q^2$-homogeneous the subspace $J_\al$ of $\uqbp$ is
  also $Q^2$-homogeneous by Relation~\eqref{eq:Q2kow}. Let $\gamma\in Q$ be
  such that $x_\gamma:=\pr _{(\gamma,-\beta)}(x) \neq 0$. Then
  $\al-\gamma\in L^\perp$ by Lemma \ref{le:adTLw}. We now proceed
  indirectly to prove the first statement. If $x_\gamma\notin J_\al$ and
  $\gamma$ is maximal with respect to $\prec$ with this property, then
  Lemma~\ref{le:tgamma} for $J=J_\al $
  implies $K_{\gamma-\beta}\in C$. Similarly, $K_{\alpha-\beta}\in C$
  by Lemma~\ref{le:tgamma} for $J=0$. Hence $\gamma-\al\in L$ and thus
  $\al=\gamma$ as $(\cdot,\cdot)$ is positive definite and
  $\al-\gamma\in L^\perp$. This however is a contradiction to
  $x_\al\in J_\al$, $x_\gamma\notin J_\al$. Clearly, $x_\gamma \in J_\al $
  implies that $\al -\gamma \in Q_+$.

  Statement (2) immediately follows from (1). To prove statement (3) assume that
  $\pr _{(\gamma ,-\beta )}(x)\not=0$ for some $\gamma \in \beta 
  +L^\perp $. By (1) one gets $\al-\beta\in L^\perp$. On the other
  hand Lemma \ref{le:tgamma} implies that $\al-\beta\in L$ and hence
  $\al=\beta$. 
\end{proof}

\subsection{Connected right coideal subalgebras}
In analogy to the terminology for coalgebras \cite[5.1.5]{b-Montg93}
we make the following definition. 

\begin{defi}
Let $C$ be a right coideal subalgebra of $U^{\ge 0}$. We say that $C$ is
{\em connected} if $C\cap U^0=k1$.
\end{defi}

Let $\rcs $ denote the set of right coideal subalgebras $C$ of $\uqbp $ such
that $C\cap U^0$ is a Hopf algebra.
For any $C\in \rcs$ we write $L(C)$ for the subgroup of $Q$ corresponding
to the Hopf subalgebra $C\cap U^0$ of $U^0$. Clearly, $C$ is connected
if and only if $L(C)=0$.

We will now show that any right coideal subalgebra $C\in\rcs$ decomposes
into the product of $T_{L(C)}$ and a connected right coideal subalgebra. To
this end define for any $C\in \rcs $ a subspace
\begin{align}
  I(C)=\mathop{\oplus}_{\beta\in Q_+} \{x\in C_{-\beta}\,|\, (\ad
  K_\al )(x)=q^{(\al ,\beta )}x\text{ for all $\al\in L(C)$}\}.
  \label{eq:IC}
\end{align}
The definition implies directly that
\begin{align}
  I(C)=\mathop{\oplus}_{\beta \in Q_+}\{x\in C_{-\beta }\,|\,\pr _{(\gamma
  ,-\beta )}(x)=0\quad \text{for all $\gamma \in Q\setminus (\beta +
  L(C)^\perp )$}\}.
  \label{eq:IC1}
\end{align}

\begin{lem}
  Let $C\in \rcs $, $\beta \in Q_+$, and let $x\in I(C)\cap C_{-\beta }$ with
  $x\not=0$.
  Then $\pr _{(\beta ,-\beta )}(x)\not=0$. If $\pr _{(\gamma ,-\beta
  )}(x)\not=0$ for some $\gamma \in Q$ then $\beta -\gamma \in
L(C)^\perp \cap Q_+$.
  \label{le:ICchar}
\end{lem}

\begin{proof}
  Let $\al $ be the unique maximal element in $Q$ with
  $\pr _{(\al ,-\beta )}(x)\not=0$, see Lemma~\ref{le:wvec}(2).
  Then $\al \in \beta +L(C)^\perp $ by Equation \eqref{eq:IC1}, and
  hence $\al =\beta $ by Lemma~\ref{le:wvec}(3). The remaining claim holds by
  Lemma~\ref{le:wvec}(1).
\end{proof}

\begin{prop} \label{pr:IC}
  Let $C\in \rcs $.
  
  (1) $I(C)$ is a connected, $\ad\,T_{L(C)}$-stable right coideal
  subalgebra of $\uqbp $.
 
  (2) The decomposition $I(C)=\mathop{\oplus} _{\beta \in Q_+}I(C)_{-\beta }$
  is an algebra grading.

  (3) The multiplication map $I(C)\ot U^0\rightarrow \uqbp$ is injective.

  (4) The multiplication map $I(C)\ot T_{L(C)}\to C$ is bijective.

  (5) Let $D$ be a connected $\ad \,T_{L(C)}$-stable right coideal
  subalgebra of $\uqbp $. If $D\subseteq C$ then $D\subseteq I(C)$. 
  If $DT_{L(C)}=C$ then $D=I(C)$.
\end{prop}

\begin{proof}
  (1) The $Q^2$-grading of $\uqbp $
  is an algebra grading, and hence $I(C)$ is a subalgebra of
  $C$. Relation~\eqref{eq:Q2kow} and Equation~\eqref{eq:IC1}
  imply that $\kow (I(C))\subseteq I(C)\ot \uqbp $.
  Let $\beta \in Q$ and let $x\in I(C)\cap C_{-\beta }\cap U^0$. If $x\not=0$
  then $\pr _{(\beta ,-\beta )}(x)\not=0$ by Lemma~\ref{le:ICchar}, and hence
  $\beta =0$. Thus $I(C)\cap U^0=\field 1$,
  and hence $I(C)$ is a connected right
  coideal subalgebra of $C$. Further, $I(C)$ is $\ad
  \,T_{L(C)}$-stable by Definition~\eqref{eq:IC}.

  (2) This is a special case of Lemma~\ref{le:Qdeg}.

  (3) Consider elements $x_\beta\in
  I(C)_{-\beta}$ and $t_{\beta}\in U^0\setminus\{0\}$ such that
  $x_\beta\neq 0$ for finitely many $\beta\in Q_+$. Choose a maximal
  $\gamma\in Q_+$ such that $x_\gamma\neq 0$. Then
  $\pr _{(\gamma,-\gamma)}(x_\gamma)\neq 0$ but
  $\pr _{(\gamma,\delta)}(x_\beta)= 0$ for all $\beta\in Q_+\setminus\{\gamma\}$,
  $\delta \in Q$ by Lemma \ref{le:ICchar}. Hence
  $\pr _{(\gamma,\delta)}(\sum_\beta x_\beta t_\beta) = \pr _{(\gamma,\delta)}
  (x_\gamma t_\gamma)\neq 0$ for a suitable $\delta\in Q$. This implies
  $\sum_\beta x_\beta t_\beta\neq 0$ which proves injectivity. 

  (4) Injectivity follows from (3). To verify surjectivity consider an
  element $x\in C_{-\beta}$ for 
  some $\beta\in Q$ and assume that $x$ is a weight vector for
  $\ad\,T_{L(C)}$. As in Lemma \ref{le:wvec} let $\alpha\in Q_+$ be
  maximal such that $\pr_{(\alpha,-\beta)}(x)\neq 0$. Lemma
  \ref{le:tgamma} implies that $K_{\alpha-\beta}\in
  T_{L(C)}$. Moreover, $x=xK_{\beta-\al }K_{\al -\beta }$ and
  $xK_{\beta -\al } \in C_{-\al }\cap I(C)$ by Equation \eqref{eq:IC}. Hence
  $x$ lies indeed in the image of the multiplication map.

  (5) Assume that $D\subseteq C$. By \eqref{eq:Cbeta} one has 
  $D_{-\beta }\subseteq C_{-\beta }$ for all $\beta \in Q$. 
  Further, $D_{-\beta }$ is $\ad \,T_{L(C)}$-stable for all
  $\beta \in Q$ since $D$ and $U^+K_{-\beta }$ are $\ad \,T_{L(C)}$-stable.
  Let $x\in D_{-\beta }$ be a weight vector for $\ad \,T_{L(C)}$.
  By Lemma \ref{le:wvec} there exists a unique maximal weight
  $\alpha\in Q$ such that $\pr_{(\al,-\beta)}\neq 0$. Since $D$ is
  connected Lemma \ref{le:tgamma} implies $\alpha=\beta$. 
  Thus $x\in I(C)$ by Equation~\eqref{eq:IC}.
  Since $\ad \,T_{L(C)}$ is diagonalizable,
  we conclude that $D\subseteq I(C)$.
  Then the last claim holds by (4).
\end{proof}

\begin{cor}\label{cor:inrcs}
  Let $C\in \rcs $ and let $L\subseteq Q$ be a subgroup.
  Assume that $C$ is $\ad \,T_L$-stable. Then $CT_L\in \rcs $.
\end{cor}

 \begin{proof}
   It follows from the assumptions that $CT_L$ is a right coideal
   subalgebra of $\uqbp $. Further, $CT_L=I(C)T_{L(C)+L}$ by Proposition
   \ref{pr:IC}(4), and hence $CT_L\cap U^0=T_{L(C)+L}$ by Proposition
   \ref{pr:IC}(3). Thus $CT_L\in\rcs $.  
\end{proof}

Define a linear map
\begin{align*}
  \varphi:\uqbp\rightarrow \uqbp, \quad \varphi(x)=\sum_{\beta\in
    Q}\pr_{(\beta,-\beta)}(x) 
\end{align*}
and observe that $\varphi$ is a homomorphism of right
$\uqbp$-comodules by Relation~\eqref{eq:Q2kow}.
For all $C\in \rcs $ define
\begin{align*}
  \bar{C}=\varphi (I(C))U^0.
\end{align*}
We will see in the following proposition that $\bar{C}$ is a right
coideal subalgebra of $\uqbp$ containing $U^0$ as in Theorem \ref{th:HS}.
For any $C\in \rcs$ let $\varphi _C:I(C)\rightarrow \uqbp$ be the
restriction of $\varphi$ 
to $I(C)$.

\begin{prop} \label{pr:varphi}
  Let $C\in \rcs $. Then the following hold.

  (1) $\varphi _C$ is an injective homomorphism of right $U^{\ge
    0}$-comodule algebras.

  (2) $\varphi _C(I(C))$ is a connected, $\ad \,U^0$-stable right
  coideal subalgebra of $U^{\ge 0}$.
  
  (3) $\bar{C}\in \rcs $, and
  the multiplication map $\varphi _C(I(C)) \ot U^0\to \bar{C}$ is
  bijective.
  
  (4) $\varphi _C(I(C))=S(U^+)\cap \bar{C}$.
\end{prop}

\begin{proof}
  (1) Lemma~\ref{le:ICchar} implies that $\varphi _C$ is injective and that
  \[ \pr _{(\beta ,-\beta )}(x)\pr _{(\gamma ,-\gamma )}(y)=
     \pr _{(\beta +\gamma ,-(\beta +\gamma ))}(xy)\quad
    \text{for all $x\in I(C)_{-\beta }$, $y\in I(C)_{-\gamma }$.}
  \]
  This and Proposition \ref{pr:IC}(2) imply that $\varphi _C$ is an algebra
  homomorphism. As $\varphi$ is a homomorphism of right
  $\uqbp$-comodules so is $\varphi_C$.

  (2) By (1) the subspace $\varphi _C(I(C))$ is a  right coideal
  subalgebra of $\uqbp $. It is connected since $\varphi(\uqbp)\cap
  U^0=\field 1$. It is $\ad \,U^0$-stable, since $\varphi
  _C(I(C)_{-\beta })$ is $\ad \,U^0$-stable for all $\beta \in Q_+$.

  Recall that the multiplication map $S(U^+)\ot U^0\to \uqbp $ is bijective.
  Thus
  (3) and (4) follow from (2) and since $\varphi _C(I(C))\subseteq S(U^+)$.
\end{proof}

Define a linear map $\psi :U^+\to S(U^+)$ by
\begin{align*}
   \psi (x_\beta )=q^{-(\beta ,\beta )/2}x_\beta K_\beta ^{-1}\quad
   \mbox{ for all $x_\beta \in U^+_\beta $, $\beta \in Q_+$.}
\end{align*}
The following technical Lemma will allow us to identify the connected
$\ad \,U^0$-stable right coideal subalgebras of $\uqbp $ in Theorem
\ref{th:connrcs}.  

\begin{lem} \label{le:psi}
  The map $\psi :U^+\to S(U^+)$ is an algebra isomorphism. Moreover,
  $\psi (U^+[w])=S(U^+)\cap U^+[w]U^0$ for all $w\in W$.
  \label{le:phi}
\end{lem}

\begin{proof}
  Since $S(U^+)$ is generated by the elements $S(E_\al )=-K_\al ^{-1}E_\al $,
  where
  $\al \in \Pi $, we conclude that $S(U^+)$ is spanned by the elements
  of the form $x_\beta K_\beta ^{-1}$, where $\beta \in Q_+$ and $x_\beta \in
  U^+_\beta $. Thus $\psi (U^+)=S(U^+)$.
  Since $U^+$ is graded by $Q_+$, the map $\psi $ is also
  injective. For any $\alpha,\beta\in Q_+$ and any  $x_\al \in U^+_\al
  $, $y_\beta \in U^+_\beta $ one verifies that
  $\psi (x_\al y_\beta )=\psi (x_\al )\psi (y_\beta )$ which implies
  that $\psi$ is an algebra homomorphism.
  
  Let $w\in W$. Since $U^+[w]$ is graded by $Q_+$,
  the last claim follows from the definition of $\psi $ and the above
  description of $S(U^+)$.
\end{proof}

\begin{thm} \label{th:connrcs}
  (1) Let $w\in W$. Then $\psi (U^+[w])$ is a connected $\ad
  \,U^0$-stable right coideal subalgebra of $\uqbp $.

  (2) Any connected $\ad \,U^0$-stable right coideal subalgebra of
  $\uqbp $ is of the form $\psi (U^+[w])$ for a unique $w\in W$.
\end{thm}

\begin{proof}
  (1) Let $C=U^+[w]U^0$. Then $C\in \rcs$ by Theorem \ref{th:HS}. Moreover,
  $L(C)=Q$ and hence $I(C)=S(U^+)\cap U^+[w]U^0$ 
  by Equation~\eqref{eq:IC} or Equation~\eqref{eq:IC1}.
  Therefore $I(C)=\psi (U^+[w])$
  by Lemma~\ref{le:psi}. Moreover, $\varphi _C=\id_{I(C)} $, and hence the
  claim holds by Proposition \ref{pr:varphi}(2).

  (2) Let now $D$ be a connected $\ad \,U^0$-stable right coideal subalgebra.
  Then $I(D)=D$ by Proposition~\ref{pr:IC}(4), and hence $\varphi _D=\id_D $
  by Equation~\eqref{eq:IC}. Moreover $\bar{D}=DU^0\in \rcs $ by
  Proposition \ref{pr:varphi}(3), and $D=S(U^+)\cap DU^0$ by
  Proposition \ref{pr:varphi}(4). By Theorem \ref{th:HS}, $DU^0=U^+[w]U^0$
  for a unique $w\in W$. Thus $D=\psi (U^+[w])$ by Lemma~\ref{le:psi}.
\end{proof}

\subsection{Right coideal subalgebras and characters}\label{sec:red-char}
By definition a character of an associative, unital $\field$-algebra
$A$ is an algebra 
homomorphism $\phi:A\rightarrow \field$ with $\phi(1)=1$. Let
$\Char(A)$ denote the set of characters of $A$. 
The following construction, though only considered for $\uqbp$ here,
works for right coideal subalgebras of an arbitrary Hopf algebra. 

Let $C$ be a right coideal subalgebra of $U^{\ge 0}$ and $\phi\in
\Char(C)$. Consider the subspace  
\begin{align*}
  C_\phi:=\{(\phi\ot \id)\kow(x)\,|\,x\in C\}=\{\phi(x_{(1)})x_{(2)}\,|\,x\in C\}
\end{align*}
of $U^{\ge 0}$. It follows from the coassociativity of the coproduct
that $C_\phi$ is a right coideal subalgebra of $\uqbp$ and that the
map 
\begin{align}\label{eq:Cphi-hom}
  C\rightarrow C_\phi,\qquad x\mapsto
  \phi(x_{(1)})x_{(2)}\quad\mbox{for all $x\in \uqbp$} 
\end{align} 
is a surjective homomorphism of right $\uqbp$-comodule algebras.

\begin{lem}\label{le:Cphi}
  Let $C$ be a connected right coideal subalgebra of $U^{\ge 0}$ and $\phi\in
  \Char(C)$. Then $C_\phi$ is connected and the map
  \eqref{eq:Cphi-hom} is an isomorphism. 
\end{lem} 

\begin{proof}
  Recall the standard $\N_0$-filtration of $\uqbp$. More explicitly,
  define a   $\Z$-linear homomorphism 
  $\hght:Q\rightarrow \Z$ by $\hght(\alpha)=1$ for all $\alpha\in
  \Pi$. Now set $\cF^n\uqbp=\oplus_{\beta\in Q,\,\hght(\beta)\le
    n}U^+_\beta U^0$. Note that $\cF$ is a filtration of $\uqbp$ both
  as an algebra and as a coalgebra. Assume that $C$ is connected. It
  suffices to show for any $x\in C$ and any $m\in \N_0$ the relation 
  \begin{align}
     x\in \cF^m\uqbp\setminus \cF^{m-1}\uqbp\, \Longrightarrow\,
     \phi(x_{(1)})x_{(2)}\notin \cF^{m-1}\uqbp. 
  \end{align}
  To this end write $x=x_m+ u$ with $x_m\in \oplus_{\beta\in
    Q,\,\hght(\beta)= m}U^+_\beta U^0$ and $u\in
  \cF^{m-1}\uqbp$. Since $\cF$ is a filtration of coalgebras and $C$
  is a connected right coideal subalgebra one obtains $\kow(x_m)-1\ot
  x_m\in \uqbp\ot \cF^{m-1}\uqbp$. Thus $\phi(x_{(1)})x_{(2)}-x_m\in
  \cF^{m-1}\uqbp$ which implies $\phi(x_{(1)})x_{(2)}\notin
  \cF^{m-1}\uqbp$ and concludes the proof.   
\end{proof}

\begin{rema}
  One can show that the map \eqref{eq:Cphi-hom} is bijective for all
  $C\in \rcs $. 
\end{rema}

Let $V\subseteq U^+$ be an $\ad \,U^0$-stable subspace and let $\phi \in \Hom
(V,\field )$ be a linear functional. We define
\begin{align}
  \supp \,\phi =\{\beta \in Q_+\,|\,\phi (x)\not=0\text{ for some
  $x\in U^+_\beta\cap V $}\}.
  \label{eq:supprho}
\end{align}
We will only use the notion $\supp\,\phi$ in the case where $V=U^+[w]$ for
some $w\in W$ and where $\phi$ is a character.
The following theorem is the first main result of this paper. Recall that
$\rcs$ denotes the set of all right coideal subalgebras $C$ of $\uqbp $
such that $C\cap U^0$ is a Hopf algebra. Recall moreover that $\vep$
denotes the counit of $\uqbp$.

\begin{thm}\label{thm:coid-class}
  (1) Let $C\in \rcs $. Then $\varphi _C(I(C))=\psi (U^+[w_C])$
  for a unique $w_C\in W$, and $\varphi _C^{-1}\psi :U^+[w_C]\to I(C)$ is
  an algebra isomorphism with $\varphi _C^{-1}\psi (U^+[w_C]\cap U^+_\beta
  )=I(C)_{-\beta }$ for all $\beta \in Q_+$.
  The map
  \begin{align}
    \phi_C :U^+[w_C]\to \field ,\quad
    \phi_C (x)=\vep (\varphi _C^{-1}\psi (x))\quad
    \text{for all $x\in U^+[w_C]$,}
    \label{eq:rho}
  \end{align}
  is a character, and
  $L(C)\subseteq (\supp \,\phi_C )^\perp $.

  (2) Let $w\in W$, $\phi\in \Char(U^+[w])$, and $L\subseteq
  (\supp \,\phi )^\perp $ a subgroup. Let 
  \[ D(w,\phi):=\psi(U^+[w])_{\phi\psi^{-1}}=\{\phi (\psi
  ^{-1}(x_{(1)}))x_{(2)}\,|\,x\in \psi (U^+[w])\} \] 
  and $C(w,\phi,L):=D(w,\phi)T_L$.
  Then $D(w,\phi)$ is a connected, $\ad\, T_L$-stable right coideal
  subalgebra of $U^{\ge 0}$, and $C(w,\phi,L)\in \rcs $. 

  (3) The map 
  \begin{align*}
    \rcs&\rightarrow \{(w,\phi ,L)\,|\,w\in W,\, \phi\in
    \Char(U^+[w]),\,L\,\mbox{is a subgroup of }\,(\supp\, \phi)^\perp
    \} \\
    C&\mapsto (w_C,\phi_C,L(C)) 
  \end{align*}
  given by (1) is a bijection with inverse map $(w,\phi,L)\mapsto
  C(w,\phi,L)$ as in (2). 
\end{thm}

\begin{proof}
  (1) By Proposition \ref{pr:varphi}(2), $\varphi _C(I(C))$ is a connected $\ad
  \,U^0$-stable right coideal subalgebra of $\uqbp $, and hence
  $\varphi _C(I(C))=\psi (U^+[w_C])$ for a unique $w_C\in W$ by
  Theorem \ref{th:connrcs}. Now $\varphi _C$ is an injective algebra
  homomorphism by Proposition \ref{pr:varphi}(1), and hence $\varphi
  _C^{-1}\psi $ is an algebra isomorphism. The compatibility of
  $\varphi^{-1}_C\psi$ with the given degrees holds since $\psi 
  (U^+_\beta )=U^+_\beta K_\beta ^{-1}$
  and since
  $\varphi _C(I(C)_{-\beta })\subseteq U^+_\beta K_\beta ^{-1}$
  for all $\beta \in Q_+$.
  Moreover, $\phi_C \in \Char (U^+[w_C])$ since
  $\vep $ and $\varphi _C^{-1}\psi $ are algebra homomorphisms.
  Finally, let $\gamma \in Q$ and $x_\gamma \in U^+[w_C]\cap U^+_\gamma $.
  Then $y:=\varphi _C^{-1}\psi (x_\gamma )\in I(C)_{-\gamma }$.
  Hence, if $\phi_C (x_\gamma )=\vep(y)=\pr _{(0,-\gamma )}(y)K_\gamma\not=0$
  then $\gamma \in L(C)^\perp $ by Lemma~\ref{le:ICchar}. Thus
  $L(C)\subseteq (\supp \,\phi_C )^\perp $.

  (2) Lemma \ref{le:Cphi} implies that $D(w,\phi)$ is a connected right
  coideal subalgebra of $U^{\ge 0}$. Moreover, $D(w,\phi)$ is $\ad
  \,T_L$-stable since the adjoint action of 
  $U^0$ commutes with $\psi $ and since $L\subseteq (\supp \,\phi
  )^\perp$. Hence $D(w,\phi)T_L\in \rcs $ by Corollary \ref{cor:inrcs}.  

  (3) Let $C\in \rcs $. For all $x\in \psi (U^+[w_C])$ one obtains 
  \[ \phi_C (\psi ^{-1}(x_{(1)}))x_{(2)}=\vep (\varphi _C^{-1}(x_{(1)}))x_{(2)}
  =\vep (\varphi _C^{-1}(x)_{(1)})\varphi _C^{-1}(x)_{(2)}=\varphi _C^{-1}(x),
  \]
  where the second equation is satisfied by Proposition \ref{pr:varphi}(1).
  Therefore
  $D(w_C,\phi_C)=\varphi _C^{-1}\psi (U^+[w_C])=I(C)$,
  and hence $D(w_C,\phi_C)T_{L(C)}=C$ by Proposition \ref{pr:IC}(4). 
  
  Conversely, consider a triple $(w,\phi ,L)$ contained in the codomain of the
  map in (3), and let $C=C(w,\phi,L)$. 
  Then $I(C)=D(w,\phi)$ by Proposition~\ref{pr:IC}(5) and $\varphi
  _C(D(w,\phi))=\psi (U^+[w])$  by definition of $\varphi _C$ and
  $D(w,\varphi)$. Thus we
  have shown $w_C=w$ and $L(C)=L$. Further, one has
  \begin{align}\label{eq:y}
    \varphi _C(\phi (\psi ^{-1}(y_{(1)}))y_{(2)})=y\quad\mbox{for all  $y\in \psi
  (U^+[w])$}
  \end{align}
  by definition of $\varphi _C$. Let $x\in U^+[w]$ 
  and let $y=\psi (x)$. Then
  \[ \phi_C(x)=\vep (\varphi _C^{-1}\psi (x))=\vep (\varphi _C^{-1}(y))
  =\vep (\phi (\psi ^{-1}(y_{(1)}))y_{(2)})=\phi (\psi ^{-1}(y))=\phi
  (x),\]
  where the first equation follows from \eqref{eq:rho} and the third
  equation holds by \eqref{eq:y}.
  Hence we have $\phi_C=\phi$ which completes the proof of (3).
\end{proof}

\section{Characters of $\Uq^+[w]$}\label{characters}
Motivated by Theorem \ref{thm:coid-class} we now turn to the
classification of characters of $U^+[w]$. Throughout this section we
fix an element $w\in W$ and a reduced expression  
\begin{align}\label{eq:w-red}
  w=s_{\alpha_1}\cdots s_{\alpha_t},\quad \al _1,\dots,\al _t\in \Pi ,
\end{align}
of $w$ in terms of simple reflections, where $t=\ell (w)$. 
As in Section~\ref{sec:qea}, for all $i=1,\dots,t$ let
$\beta_i=s_{\alpha_1}\cdots s_{\alpha_{i-1}}\alpha_i$. The set
$\{\beta_i\,|\,i=1,\dots, t\}$ coincides with 
\begin{align*}
  \Phi_w^+:=\{\beta\in \Phi^+\,|\,w^{-1}\beta\in \Phi^-\}.
\end{align*}
Recall that $\Uq^+[w]$ is the subalgebra of $\Uq^+$ generated by the
root vectors $E_{\beta_i}:=T_{\alpha_1}\cdots T_{\alpha_{i-1}}
E_{\alpha_i}$ for $i=1,\dots,t$. In the following, if we write
$E_\beta$ for some $\beta\in \Phi_w^+$ we always mean the root vector
corresponding to the fixed reduced expression \eqref{eq:w-red} for the
element $w$. If we write $\beta_i$ with a lower index $i=1,\dots,t$
then we always refer to the specific ordering of $\wurz_w^+$ from
above. 
\subsection{Orthogonality}
Let $\phi:\Uq^+[w]\rightarrow \qfield$ be a character. In the
following lemma we show that roots corresponding to root vectors
on which $\phi$ does not vanish, are mutually orthogonal.
By the result of Levendorski{\u \i} and Soibelman \cite[Proposition
5.5.2]{a-LevSoi91}, \cite[p.~164]{b-Jantzen96} 
one has for $i<j$ the commutation relation
\begin{align}\label{LS}
  E_{\beta_i}E_{\beta_j} - q^{(\beta_i,\beta_j)}
  E_{\beta_j}E_{\beta_i}=\sum_{\makebox[2cm]{\tiny $(a_{i+1},\dots,a_{j-1})\in
  \N_0^{j-i-1}$}} m_{(a_{i+1},\dots,a_{j-1})} 
  E_{\beta_{i+1}}^{a_{i+1}}\cdots  E_{\beta_{j-1}}^{a_{j-1}}.
\end{align}
for some coefficients  $m_{(a_{i+1},\dots,a_{j-1})} \in
\qfield$. These commutation relations are the main ingredient to
obtain the following result. 
\begin{lem}\label{orthogonalLem}
  Let $\phi:\Uq^+[w]\rightarrow \qfield$ be a character and
  $(\beta_i,\beta_j)\neq 0$ for some $i\neq j$. Then
  $\phi(E_{\beta_i})=0$ or $\phi(E_{\beta_j})=0$. 
\end{lem}
\begin{proof}
  Consider Equation \eqref{LS} for $i<j$. If the right hand side is
  zero, then the claim of the lemma holds. We now perform an indirect
  proof. Assume that $j-i>0$ is minimal such that
  $(\beta_i,\beta_j)\neq 0$ and such that $\phi(E_{\beta_i})\neq 0$
  and $\phi(E_{\beta_j})\neq 0$. Then the right hand side of
  Equation~\eqref{LS} is nonzero and contains a monomial
  $E_{\beta_{i+1}}^{a_{i+1}}\cdots  E_{\beta_{j-1}}^{a_{j-1}}$ such
  that $\phi(E_{\beta_m})\neq 0$ whenever $i<m<j$ and $a_m\neq
  0$. Thus by Lemma~\ref{le:nonorth}
  we can find an $n$ with $i<n<j$ such that
  $\phi(E_{\beta_n})\neq 0$ and one of $(\beta_i,\beta_n)$,
  $(\beta _j,\beta _n)$ is nonzero.
  This is a contradiction to the minimality of $j-i>0$.
\end{proof}
For any $\phi\in \Char(U^+[w])$ define
\begin{align} \label{eq:Deltaphi}
    \wurz_w^+(\phi):=\{\beta\in \wurz_w^+\,|\,\phi(E_\beta)\neq 0\}.
\end{align} 
We will see in Remark \ref{re:independent} that $\wurz_w^+(\phi)$ is
independent of the chosen reduced expression for $w$.
\begin{rema}\label{re:lin-indep}
By Lemma \ref{orthogonalLem} for any  $\phi\in \Char(U^+[w])$ the set
$\wurz_w^+(\phi)$ consists of pairwise orthogonal roots. In
particular, $\wurz_w^+(\phi)$ is linearly independent and contains at
most $\rank(\gfrak)$ elements. 
\end{rema}
\subsection{Polynomial $H$-prime ideals}\label{sec:H-poly}
The algebra $\Uq^+[w]$ is $Q$-graded via the direct sum decomposition 
$\Uq^+[w]=\oplus _{\al \in Q}(\Uq^+[w]\cap U^+_\al )$.
Equivalently, $U^+[w]$ has a rational action of the torus
$H:=(\qfield^\times)^{\rank (\gfrak )}$ by $\qfield $-algebra
automorphisms, where $\field^\times=\field\setminus\{0\}$,  
(cf. \cite[II.2.11]{b-BG02}). This will allow us to follow
\cite[Chapter II]{b-BG02} and apply the theory of prime ideals invariant
under a torus action to determine $\Char(U^+[w])$.
Recall that a proper ideal $P$ of an associative $\field$-algebra $A$
is called prime if any two ideals $I$, $J$ not contained in $P$
satisfy $IJ \not\subseteq  P$. Equivalently, $P$ is prime if for all 
$a,b\in A\setminus P$ there exists $c\in A$ with $acb\notin P$. In
particular, if $P$ is prime and $A/P$ is commutative, then $a,b\in
A\setminus P$ implies that $ab\notin P$. 
If $A$ admits an action of the torus $H$ by automorphisms, then an ideal $P$ of
$A$ is called $H$-prime, see \cite[II.1.9]{b-BG02}, 
if any two ideals $I$, $J$ stable under the action of $H$ and not
contained in $P$ satisfy $IJ \not\subseteq P$. If, moreover, $A$ is
noetherian and $H$ acts rationally, then $H$-prime ideals are
prime \cite[II.2.9]{b-BG02}. All these assumptions are fulfilled for
the action of $H=(\field^\times)^{\rank \gfrak}$ on $U^+[w]$.

Let $\phi:\Uq^+[w]\rightarrow \qfield$ be a character. Then
$\ker(\phi)$ is an ideal in $\Uq^+[w]$ which is one-codimensional and hence
prime. The intersection  
\begin{align}\label{PM}
  P_\phi:=\bigcap_{h\in H} h\ker(\phi)
\end{align}
is an $H$-prime ideal of  $\Uq^+[w]$. Indeed, let $I$ and $J$ be
$H$-stable (or equivalently $Q$-graded) ideals
such that $IJ\subseteq P_\phi\subseteq \ker(\phi)$. Then, say,
$I\subseteq  \ker(\phi)$ because $\ker(\phi)$ is prime and
hence $I\subseteq P_\phi $ since $I$ is $H$-stable. Following
\cite[II.1.9]{b-BG02}  the set of $H$-prime ideals in $U^+[w]$ will be
denoted by $H-\spec(U^+[w])$.  
The $H$-prime ideals of $\Uq^+[w]$ obtained in the above way
via characters have very special properties as illustrated by the
following lemma. 
\begin{prop}\label{prop:poly}
  Let $\phi:\Uq^+[w]\rightarrow \qfield$ be a character and $P_\phi$
  the corresponding 
  $H$-prime ideal defined by Equation~\eqref{PM}. Then 
  the quotient $\Uq^+[w]/ P_\phi$ is a commutative polynomial
  ring in the root vectors $E_{\beta}$ with $\beta\in \wurz_w^+(\phi)$.
\end{prop}
\begin{proof}
  The quotient algebra  $\Uq^+[w]/ P_\phi$ is generated by the root
  vectors $E_{\beta}$ with $\beta\in \wurz_w^+(\phi)$. 
  Assume now that $\phi(E_{\beta_i})\neq 0$ and $\phi(E_{\beta_j})\neq 0$ for
  some $i< j$. We know from Lemma \ref{orthogonalLem} that
  $(\beta_i,\beta_j)=0$. Hence by Equation~\eqref{LS} the commutator of
  $E_{\beta_i}$ and $E_{\beta_j}$ is a polynomial in the
  $E_{\beta_m}$ of weight $\beta_i+\beta_j$. As
  $(\beta_i,\beta_i+\beta_j)=(\beta_i,\beta_i)\neq 0$ each
  monomial in this polynomial contains a root vector $E_{\beta_m}$
  with $m\neq i$ and $(\beta_i,\beta_m)\neq 0$. By Lemma
  \ref{orthogonalLem} this implies $\phi(E_{\beta_m})=0$ and hence
  $E_{\beta_m}\in P_\phi$. Thus, the right hand side of Equation~\eqref{LS} is
  contained in $P_\phi$ and the commutator $[E_{\beta_i},E_{\beta_j}]$
  is zero in $\Uq^+[w]/ P_\phi$. This proves that $U^+[w]/P_\phi$ is commutative.

  It remains to show that in $\Uq^+[w]/ P_\phi$ there are no relations
  between the $E_{\beta}$ with $\beta\in \wurz^+_w(\phi)$ apart
  from commutativity. To this end assume that $P_\phi$ contains a nonzero
  commutative polynomial in the $E_{\beta}$ with $\beta\in
  \wurz_w^+(\phi)$. By linear independence of the set
  $\wurz_w^+(\phi)$ observed in Remark \ref{re:lin-indep} one
  can use the torus action (or equivalently the $Q$-grading) to show that 
  $P_\phi$ already contains a monomial in the $E_{\beta}$ with $\beta\in
  \wurz_w^+(\phi)$. As noted above $P_\phi$ is a prime ideal and
  $\Uq^+[w]/P_\phi$ is commutative. Hence 
  one gets the contradiction $E_{\beta}\in P_\phi$ for some
  $\beta\in \wurz_w^+(\phi)$. 
\end{proof}
\begin{rema}\label{re:independent}
  By Proposition \ref{prop:poly} the set $\wurz_w^+(\phi)$ defined in
  \eqref{eq:Deltaphi} is uniquely determined by the $H$-prime ideal
  $P_\phi$. Hence $\wurz_w^+(\phi)$ only depends on the character $\phi$ and
  not on the chosen reduced expression for $w$.  
\end{rema}
Proposition \ref{prop:poly} motivates the following definition.
\begin{defi}\label{defi:poly}
  An $H$-prime ideal $P$ of $\Uq^+[w]$ is called {\em polynomial} if
  $\Uq^+[w]/P$  is a commutative polynomial ring in the $E_{\beta}$,
  $\beta\in \Phi^+_w$, 
  which are not contained in $P$. The set of polynomial $H$-prime
  ideals in  $\Uq^+[w]$ is denoted by $H-\spec^{poly}(U^+[w])$. 
\end{defi}
\begin{lem}\label{lem:poly-gen}
  Let $P\in H-\spec^{poly}(U^+[w])$. The ideal $P$ is generated by the set
  $\{E_\beta\,|\,\beta\in \wurz^+_w, E_\beta\in P\}$. 
\end{lem}
\begin{proof}
  This follows from the PBW-Theorem for $U^+[w]$, cf.~\cite[8.24]{b-Jantzen96}.
\end{proof}
To summarize results we collect some equivalent characterizations of polynomial
$H$-prime ideals of $\Uq^+[w]$.
\begin{prop}\label{prop:poly-primes}
  For any $H$-prime ideal $P$ of $\Uq^+[w]$ the following are equivalent:
  \begin{enumerate}
    \item The $H$-prime ideal $P$ is polynomial.
    \item The algebra $\Uq^+[w]/P$ is commutative.
    \item There exists $\phi\in\Char(\Uq^+[w])$ such that the prime
      ideal $P$ is generated by the set
      $\{E_{\beta}\,|\,\beta\in \wurz ^+_w\setminus \wurz^+_w(\phi)\}$. 
  \end{enumerate}
\end{prop}
\begin{proof}
  The equivalence of properties (1) and (3) is immediate from Lemma
  \ref{lem:poly-gen} and Proposition \ref{prop:poly}. It 
  remains to show that (2) implies (1). To this end assume that $P$ is
  an $H$-prime ideal such that $U^+[w]/P$ is commutative and define
  \begin{align*}
    J=\{i\in \{1,\dots,t\}\,|\,E_{\beta_i}\notin P\}.
  \end{align*} 
  We claim that $(\beta_i,\beta_j)=0$ for all distinct $i,j\in
  J$. Indeed, otherwise choose elements $i<j$ in $J$ with
  $j-i$ minimal such that $(\beta_i,\beta_j)\neq 0$ and
  $E_{\beta_i},E_{\beta_j}\notin P$. This implies for any $m$
  with $i<m<j$ that $E_{\beta_m}\in P$ if
  $(\beta_i,\beta_m)\neq 0$
  or $(\beta_m,\beta_j)\neq 0$.
  Hence by Lemma~\ref{le:nonorth}
  the right hand side of  Equation~\eqref{LS} belongs to $P$
  for the chosen $i,j$. The assumption
  $(\beta_i,\beta_j)\neq 0$ now implies that
  $E_{\beta_i}E_{\beta_j}\in P$.
  But $U^+[w]/P$ is commutative and $P$ is prime,
  which gives the contradiction to
  $E_{\beta_i},E_{\beta_j}\notin P$. This proves that indeed
  $(\beta_i,\beta_j)=0$ for all distinct $i,j\in J$. 

  Assume now that a commutative polynomial $f$
  in the $E_{\beta_i}$, $i\in J$, belongs to $P$.
  Using the $H$-stability of $P$ and the
  orthogonality one may assume that $f$ is a monomial. Yet this yields
  again a contradiction to the fact that $P$ is prime
  and $U^+[w]/P$ is commutative.
\end{proof}
\subsection{Stratification of $\Char(U^+[w])$}\label{sec:H-strat}
By Proposition \ref{prop:poly} one has a surjective map
\begin{align*}
  \pr_w: \Char(\Uq^+[w])\rightarrow \mbox{$H-\spec^{poly}(U^+[w])$}, \quad
  \phi\mapsto \pr_w(\phi):=P_\phi. 
\end{align*}
For each $P\in \mbox{$H-\spec^{poly}(U^+[w])$}$ the preimage
$\pr_w^{-1}(P)$ is isomorphic to the spectrum of a Laurent polynomial
ring in as many variables as there are elements in
$\{\beta\in \Phi^+_w\,|\,E_\beta\notin P\}$. Therefore the set $\Char(U^+[w])$ is
a disjoint union
\begin{align}\label{char-strat}
  \Char(U^+[w])=\bigsqcup_{P\in  H-\spec^{poly}(U^+[w])}\pr_w^{-1}(P)
\end{align}
of spectra of Laurent polynomial rings. The decomposition
\eqref{char-strat} is the $H$-stratifica-tion of $\Char(\Uq^+[w])$ in
the sense of \cite[II.2.1]{b-BG02}. In order to classify characters of
$\Uq^+[w]$ it hence remains to determine all polynomial $H$-prime
ideals $P$ of $U^+[w]$ and for each of them the set 
$\{\beta\in \Phi^+_w\,|\,E_\beta\notin P\}$.

Let $\Theta\subseteq \Phi^+$ be a set of pairwise orthogonal
roots. Then the reflections $s_\beta$ and $s_\gamma$ commute for any
$\beta ,\gamma \in \Theta$. Hence we may write $ \Pi_{\beta\in
  \Theta}s_\beta$ to denote the product of all reflections
corresponding to roots in $\Theta$. Let 
\begin{align}\label{eq:w-Theta}
  w_\Theta=(\Pi_{\beta\in \Theta}s_\beta)w
\end{align}
and define
\begin{align}\label{eq:Tw-def}
  T^w=\{\Theta\subseteq \wurz^+_w\,|\,&\mbox{$\Theta$ consists of
    pairwise orthogonal roots,}\\
    &\mbox{and }\, \ell(w_\Theta)=\ell(w)-|\Theta|\}.\nonumber
\end{align}
For any $\Theta\in T^w$ and any $\beta\in \Phi^+_w$ one has the implication
\begin{align}\label{eq:inTheta}
  \beta\in \Theta \,\Longrightarrow \, w_\Theta^{-1}\beta\in \Phi^+.
\end{align}
In explicit examples of small rank the set $T^w$ is not hard to
determine. Let $J\subseteq \{1,\dots,t\}$ be a subset such that the
elements of $\Theta:=\{\beta_i\,|\,i\in J\}$ are pairwise
orthogonal. Then the set 
$\Theta$ belongs to $T^w$ if and only if one obtains a reduced
expression by deleting all simple reflections $s_{\alpha_i}$ for $i\in
J$ from the expression \eqref{eq:w-red} for $w$. Moreover, by the
following lemma, all sets in $T^w$ can be built by an inductive
procedure from smaller sets in $T^w$. 
\begin{lem}\label{lem:w-sub-red}
   Let $w\in W$ and $\Theta\in T^w$. If $\Theta'\subseteq \Theta$ then
   $\Theta'\in T^w$. 
\end{lem}
\begin{proof}
  Let $\Theta'' \subset \Theta$ and $\beta\in \Theta\setminus \Theta''$.
  For any weight $\lambda$ the orthogonality assumption implies that
  $(\beta, w_{\Theta''}\lambda)=(\beta,w\lambda)$. Hence, if $\lambda$
  is regular dominant, that is $(\alpha,\lambda)>0$ for all $\alpha\in
  \Pi$, then $(\beta, w_{\Theta''}\lambda)<0$ because $w^{-1}\beta<0$.  
  By Lemma \ref{le:weyl-bruhat} this implies
  $s_\beta  w_{\Theta''} < w_{\Theta''}$. Repeated
  application of this argument leads to
  \begin{align*}
    \Theta'\subsetneqq \Theta''\subseteq \Theta\quad \Longrightarrow
    \quad w_\Theta\le w_{\Theta''}<w_{\Theta'}. 
  \end{align*}
  Together with $l(w_\Theta)=l(w)-|\Theta|$ this implies
  $l(w_{\Theta'})=l(w)-|\Theta'|$ 
  for any $\Theta'\subseteq \Theta$. Hence any subset
  $\Theta'\subseteq \Theta$ belongs to $T^w$.
\end{proof}
The set $T^w$ will provide us with the desired parametrization of
{$H-\spec^{poly}(U^+[w])$. As a first step we associate a polynomial
  $H$-prime ideal to each element in $T^w$.  
\begin{prop}\label{prop:w-reducing}
  Let $\Theta\in T^w$ and let $P_\Theta$ denote the two-sided ideal of
  $U^+[w]$ generated by all $E_\beta$ with $\beta\in
  \wurz^+_w\setminus \Theta$. Then  
  $P_\Theta$ is a polynomial $H$-prime ideal of
  $\Uq^+[w]$, and $E_{\beta}\notin P_\Theta$ for all $\beta\in \Theta$. 
\end{prop}
\begin{proof}
  We may assume that $\gfrak$ is simple.
  By \cite[8.24]{b-Jantzen96} the algebra $\Uq^+[w]$ can be given in
  terms of generators $E_{\beta_i}$, $i=1,\dots,t$, and relations
  \eqref{LS}. It suffices to show that for any $m,n\in \{1,\dots, t\}$
  with $m<n$ there is no family $\{a_j\in \N_0\}_{\beta_j\in \Theta}$
  such that 
  \begin{align}\label{not-pos}
    \beta_m+\beta_n=\sum_{\beta_j\in \Theta, m<j<n} a_j \beta_j.
  \end{align}  
  Indeed, if no such family exists, then any nontrivial monomial on
  the right hand side of Equation~\eqref{LS} contains a factor in
  $P_\Theta$. 
  
  If $\gfrak$ is of rank 2 then $T^w$ consists only of subsets of
  $\{\beta_1,\beta_t\}$. In this case the right hand side
  of Equation \eqref{not-pos} always vanishes. Hence we may assume that
  $\rk (\gfrak )\ge 3$. 

  We prove the impossibility of \eqref{not-pos} indirectly. Assume we
  have found a family $\{a_j\in \N_0\}_{\beta_j\in \Theta}$ such that
  \eqref{not-pos} holds. By Lemma~\ref{lem:w-sub-red} all subsets of
  $\Theta$ are also contained in $T^w$. Thus, leaving out certain
  elements of $\Theta$ we may 
  assume that all $a_j$ are strictly positive and that $\beta_m, \beta_n\notin
  \Theta$. Moreover, shortening $w$ from both sides if necessary, we may
  assume that $m=1$ and $n=t$. Hence we have
  \begin{align}\label{not-pos2}
     \alpha_1+\beta_t=\sum_{\beta_j\in \Theta} a_j \beta_j.
  \end{align} 
  Now let $\beta'_t=-w_\Theta \alpha_t \in \wurz ^+$.
  The assumptions $\Theta\in T^w$ and $\alpha _1 \notin
  \Theta$ imply that  
  \begin{align}\label{negative}
    w_\Theta^{-1} \alpha_1\in \wurz^-, \qquad
    w_\Theta^{-1}\beta_t'=-\alpha_t\in \wurz^-.
  \end{align}
  Using the definition \eqref{eq:w-Theta} of $w_\Theta$ and
  the orthogonality of the $\beta_j\in \Theta$ one calculates 
    \begin{align*}
     \beta_t'
     &=\beta_t-2 \sum_{\beta_j\in \Theta}
     \frac{(\beta_j,\beta_t)}{(\beta_j,\beta_j)} \beta_j.
  \end{align*}
  Hence one obtains that
  \begin{align}\label{nearly-contra}
    \alpha_1+\beta_t'=\sum_{\beta_j\in \Theta}\Big( a_j
    -2\frac{(\beta_j,\beta_t)}{(\beta_j,\beta_j)}\Big) \beta_j.
  \end{align}
  Apply $w_\Theta ^{-1}$ to this equation. By \eqref{negative} and since
  $w_\Theta ^{-1}\beta _j\in \wurz ^+$ for all $\beta _j\in \Theta $ by
  \eqref{eq:inTheta}, we conclude that there exists $\beta _i\in \Theta $ such
  that
  \begin{align}
    0<a_i<\frac{2(\beta _i,\beta _t)}{(\beta _i,\beta _i)}.
    \label{eq:aj}
  \end{align}
  Hence $\beta _t$ is long by Lemma~\ref{le:rank2roots}.
  Replacing $w$ by $w^{-1}$ and $\Theta $ by $-w^{-1}(\Theta )$ and applying
  the same procedure one obtains that $\al _1$ has to be long.
  Again by Lemma~\ref{le:rank2roots} the root $\beta _i\in \Theta $
  satisfying \eqref{eq:aj} has to be short.
  Then $(\alpha _1,\beta _i),(\beta _i,\beta _t)\in \{2,-2\}$, since
  $\gfrak $ is irreducible and not of type $G_2$. Further, 
  $(\alpha_1,\beta_i)+(\beta_t,\beta_i)=2a_i>0$ by
  Equation~\eqref{not-pos2}, and hence 
  $(\alpha_1,\beta_i)=(\beta_t,\beta_i)=2$. This is a contradiction to Relation
  \eqref{eq:aj}.
\end{proof}
\begin{rema}
  We will see in Corollary \ref{cor:summary} that for any $\Theta\in T^w$ the
  polynomial $H$-prime ideal $P_\Theta$ obtained by Proposition
  \ref{prop:w-reducing} is independent of the chosen reduced
  expression for $w$.
\end{rema}
\subsection{$H$-prime ideals with noncommutative quotients}
The next aim is to show that Proposition \ref{prop:w-reducing} already
produces as many distinct polynomial $H$-prime ideals as we can
possibly find. To this end we resort to the description of
$H-\spec(U^+[w])$ recently given by Yakimov \cite{a-Yakimov09p} based
on results by Gorelik \cite{a-Gorelik00}. 
In order to refer to these papers without much rewriting,
we work with $U^-[w]$ instead of $U^+[w]$ for now. Recall from
\cite[8.24]{b-Jantzen96} that by definition $U^-[w]=\omega(U^+[w])$
where $\omega$ denotes the algebra automorphism of $\Uq$ defined by
\begin{align}\label{eq:om-def}
  \omega(E_\alpha)&=F_\al,& \omega(F_\al)&=E_\al, &
  \omega(K_\al)=K_\al^{-1} \quad\mbox{for all $\al\in \Pi$}.
\end{align}

Let $\weights$ denote the weight lattice and $\weights^+$ the set of
dominant weights. 
In this subsection we mostly use Gorelik's and Joseph's notation following
\cite{a-Gorelik00}, \cite{b-Joseph}. For any dominant integral weight
$\lambda\in \weights^+$ let $V(\lambda)$ denote the corresponding simple
$U$-module and let  $v_\lambda\in V(\lambda)$ be a highest weight
vector. We write $V(\lambda)^\ast$ to denote the dual space of $V(\lambda)$
which is canonically a right $\uqg$-module. Let $\cqg$ be the Hopf
algebra generated by all matrix coefficients $c_{f,v}^\lambda$ for
$v\in V(\lambda)$, $f\in V(\lambda)^\ast$ of the representations
$V(\lambda)$, $\lambda\in \weights^+$. We define  
\begin{align*}
  R^+=\bigoplus_{\lambda\in \weights^+} V(\lambda)^\ast
\end{align*}
which is an algebra with the Cartan multiplication. We can consider
$R^+$ as a subalgebra of $\cqg$ if we identify $f\in V(\lambda)^\ast$
with the matrix coefficient $c^\lambda_{f,v_\lambda}$. Let
$\xi_{w\lambda}\in V(\lambda)^\ast $ be a functional which only   
lives on the weight space of weight $w\lambda$ in $V(\lambda)$. To
shorten notation we write $c^\lambda_{\xi_{w\lambda},v_\lambda}$
simply as $c^\lambda_w$. The elements $c^\lambda_w$ are defined up to
scalars. By \cite[9.1.10]{b-Joseph} these scalars can be chosen such
that $c^\mu_w c^\nu_w=c^{\mu+\nu}_w$ for any $\mu,\nu\in \weights^+$ and
$c_w:=\{c_w^\lambda\,|\, \lambda \in \weights^+\}$ becomes an Ore set in
$R^+$. Recall that we have fixed $w$ for all of this section, however,
we can define $c_y^\lambda$ and $c_y$ analogously for any $y\in W$.
Consider the localized algebra $R^w:=R[c_w^{-1}]$. This algebra is
$\weights$-graded with $\deg((c_w^\mu)^{-1}c^\nu_{f,v_\nu })=\nu-\mu$ for all
$\mu ,\nu \in \weights^+$, $f\in V(\nu )^*$. Let 
$R^w_0$ denote the algebra of elements of degree zero in $R^w$. Using
the abbreviation $c_w^{-\lambda}:=(c_w^\lambda)^{-1}$ one
may write
\begin{align*}
  R_0^w=\sum_{\lambda\in \weights^+}
  c_w^{-\lambda}V(\lambda)^*=\lim_{\overrightarrow{\lambda\in
      \weights^+}}c_w^{-\lambda} V(\lambda)^\ast. 
\end{align*}
In particular, any element of $R^w_0$ can be written in the
form $c_w^{-\lambda} c^\lambda_{f,v_\lambda }$ for some $\lambda \in
\weights^+$ and $f\in V(\lambda)^*$. 

Following \cite[5.2.3, 6.1.2]{a-Gorelik00}, for any $y \in W$ define
ideals $Q(y)^\pm $ of $R^+$ and $Q(y)_w^{\pm}$ of $R_0^w$ by 
\begin{align*}
  Q(y)^\pm = &\sum _{\lambda \in \weights^+}\{
  c^\lambda_{\xi,v_\lambda}\,|\, \xi
  \in V(\lambda)^\ast, \, \xi\perp \Uq^\pm v_{y\lambda}\},\\
  Q(y)^\pm_w = &\lim _{\overrightarrow{\lambda\in \weights^+}} 
  \{c^{-\lambda}_w c^\lambda_{\xi,v_\lambda}\,|\, \xi
  \in V(\lambda)^\ast, \, \xi\perp \Uq^\pm v_{y\lambda}\},
\end{align*}
where $v_{y\lambda}\in V(\lambda)$ is a weight vector of weight
$y\lambda$. (For the fact that $Q(y)^\pm $ and $Q(y)^\pm_w$ are ideals confer to
Gorelik's reference to \cite[10.1.8]{b-Joseph}. Consult also
\cite[3.1]{a-Yakimov09p}.) By definition, for any $y,y'\in W$ with
$y\le y'$ one obtains the inclusions $Q(y')_w^+\subseteq Q(y)_w^+$ and 
$Q(y)_w^-\subseteq Q(y')_w^-$.

For our purposes it is sufficient to consider for any $y\le w$ the
subspaces
\begin{align}\label{eq:Qyww}
  Q(y,w)_w:=&Q(y)_w^- + Q(w)^+_w \subset R^w_0,\\
  Q(y,w):=& \{a\in R^+\,|\,\exists \lambda \in \weights^+ :
  c^\lambda _wa\in Q(y)^-+Q(w)^+\}\subset R^+.
\end{align}
For any $y\le w$ the subspace $Q(y,w)_w$ is an $H$-stable prime ideal
of $R^w_0$ by 
\cite[6.6]{a-Gorelik00}, and $Q(y,w)$ is an $H$-stable prime ideal
of $R^+$ by \cite[6.7]{a-Gorelik00}.
\begin{lem} \cite[Lemma 6.10]{a-Gorelik00} \label{le:GorLem6.10}
  Let $y,y'\in W$ with $y\le y'\le w$. Then $Q(y,w)\cap c_{y'}=\emptyset $.
\end{lem}
\begin{cor} \label{co:not0Lem}
  Let $\lambda \in \weights^+$ and $y,y'\in W$ with $y\le y'\le w$.
  Then $c_w^{-\lambda} c_{y'}^\lambda$ is not contained in $Q(y,w)_w$.
\end{cor}
\begin{proof}
  Assume that $c_w^{-\lambda} c_{y'}^\lambda \in Q(y,w)_w$ for some
  $y\le y'\le w$. Then there exists $\mu\in \weights^+$
  such that
  \begin{align*}
    c_w^{-\lambda} c_{y'}^\lambda = c_w^{-\mu} c_{\xi, v_\mu}^\mu +
    c_w^{-\mu} c_{\zeta, v_\mu}^\mu 
  \end{align*}
  for some $\xi\in (\Uq^- v_{y\mu})^\perp$ and $\zeta\in (\Uq^+
  v_{w\mu})^\perp$. Multiplying from the left by $c_w^{\lambda+\mu}$
  and setting $\eta=\lambda+\mu$ one obtains
  (using that $Q(y)^-$ and $Q(w)^+$ are ideals)
  that
  \begin{align*}
    c_w^\mu c_{y'}^\lambda =
    c_{\xi', v_\eta}^\eta +  c_{\zeta', v_\eta}^\eta 
  \end{align*}
  for some $\xi'\in (\Uq^- v_{y\eta})^\perp$ and $\zeta'\in (\Uq^+
  v_{w\eta})^\perp$. By definition of $Q(y,w)$ this means that
  $c^\lambda_{y'}\in Q(y,w)$ which is
  a contradiction to Lemma~\ref{le:GorLem6.10}.
\end{proof}
The following proposition is the main technical tool to show that
certain $H$-prime ideals of $U^+[w]$ are not polynomial. 
\begin{prop}\label{prop:not-com}
  Let $\beta,\gamma\in \wurz^+$. Define $y'=s_\gamma w$ and $y=s_\beta
  s_\gamma w$ and assume that    $y<y'<w$. Then in 
  $R^w_0/Q(y,w)_w$ the relation
  \begin{align}\label{eq:qCcommute}
     (c^{-\lambda}_w c^\lambda_{y}) (c^{-\lambda}_w c^\lambda_{y'})
    &= q^{C(\lambda)} (c^{-\lambda}_w
    c^\lambda_{y'}) (c^{-\lambda}_{w} c^{\lambda}_{y})
  \end{align}
  holds for all $\lambda\in \weights^+$, where
  $C(\lambda)= (w\lambda-y'\lambda, y\lambda- y'\lambda)$.
  In particular, if  $(\beta,\gamma)\neq 0$ then
  $R^w_0/Q(y,w)_w$ is not commutative.  
\end{prop}
\begin{proof}
  Following  \cite[Section 4]{a-Gorelik00} we use commutation
  relations in $\cqg$ which appeared as standard tools for instance in
  \cite[9.1]{b-Joseph}. 
  For any $\eta \in \weights$ define $J^+_\lambda(\eta)_w$
  (resp.~$J^-_\lambda(\eta)_w$) to be the left
  ideal of $R^w_0$ generated by $c^{-\lambda}_w c^\lambda_{f_{\eta'},v_\lambda}$ with
  $\eta'<\eta$ (resp. $\eta'>\eta$), where $f_{\eta'}\in V(\lambda) ^*$ only
  lives on the weight space of  weight $\eta'$ in $V(\lambda)$. By
  definition one has for all $w'\in W$ and $\lambda\in \weights^+$ the
  inclusions 
  \begin{align*}
    J_\lambda^+(w'\lambda)_w \subseteq Q(w')^+_w,\qquad   J_\lambda^-(w'\lambda)_w
    \subseteq Q(w')^-_w
  \end{align*}
  and hence definition \eqref{eq:Qyww} gives
  \begin{align}\label{eq:inclusion}
      J_\lambda^+(w \lambda)_w + J_\lambda^-(y\lambda)_w \subseteq Q(y,w)_w.
  \end{align}  
  By \cite[Lemma 4.2(i)]{a-Gorelik00} one has
  \begin{align}\label{eq:Gor1}
    c_w^{-2\lambda} c_{y'}^\lambda c^\lambda_w =
    q^{(\lambda,\lambda)-(w\lambda,y'\lambda)} c^{-\lambda}_w
    c^\lambda_{y'} \quad \mbox{mod }   J_\lambda^+(w \lambda)_w
  \end{align}
  and by \cite[Lemma 4.4(iv)]{a-Gorelik00} one has
  \begin{align}\label{eq:Gor2}
     (c^{-\lambda}_w c^\lambda_{y}) (c^{-\lambda}_w
     c^\lambda_{y'})= q^{-(y\lambda,y'\lambda-w\lambda)}
       (c^{-2\lambda}_w c^\lambda_{y'} c^\lambda_w )(c^{-\lambda}_w
     c^\lambda_{y}) \quad \mbox{mod } J_\lambda^-(y\lambda)_w.
   \end{align}
   The inclusion \eqref{eq:inclusion} together with
   Equations \eqref{eq:Gor1}
   and \eqref{eq:Gor2} imply the relation
    \begin{align*}
      (c^{-\lambda}_w c^\lambda_{y}) (c^{-\lambda}_w c^\lambda_{y'})
      &= q^{C(\lambda)}
      (c^{-\lambda}_w c^\lambda_{y'}) (c^{-\lambda}_{w}
      c^{\lambda}_{y}) \quad \mbox{mod } Q(y,w)_w
    \end{align*}
   with 
   $C(\lambda)=
   -(y\lambda,y'\lambda-w\lambda)
   +(\lambda,\lambda)-(w\lambda,y'\lambda)
   =(w\lambda-y'\lambda, y\lambda- y'\lambda)$. 
   This proves the commutation relation \eqref{eq:qCcommute}.

  To prove the second statement recall from Corollary \ref{co:not0Lem} that
  both $c^{-\lambda}_w c^\lambda_y$  and $c^{-\lambda}_w
  c^\lambda_{y'}$ are nonzero elements in $R^w_0/Q(y,w)_w$. We now
  assume that $\lambda$ is a regular weight, i.e.~that the Weyl group
  acts faithfully on $\lambda$. Then $w\lambda-y'\lambda$ is a nonzero
  multiple of $\gamma$, and similarly 
  $y'\lambda-y\lambda$ is a nonzero multiple of $\beta$. Hence $C(\lambda)$
  is a nonzero multiple of $(\beta,\gamma)$
  and therefore nonzero by assumption. Since
  $Q(y,w)_w$ is prime, we conclude that $R^w_0/Q(y,w)_w$ is
  not commutative. 
\end{proof}
\subsection{Classification of characters}\label{sec:char-class}
By \cite[Theorem 3.8]{a-Yakimov09p} there exists an order preserving bijection
between the poset of $H$-prime ideals of $U^-[w]$ ordered by inclusion
and the set 
\begin{align}
  W^{\le w}:=\{y\in W\,|\, y\le w\}
\end{align}
with the Bruhat order. Let $P^-(y)$ denote the $H$-prime ideal of  $U^-[w]$
corresponding to $y\in W^{\le w}$. More explicitly,
Yakimov constructs a surjective algebra homomorphism
\begin{align}
  \phi_w:R^w_0\rightarrow U^-[w]
\end{align}
with kernel $Q(w)_w^+$ which is moreover compatible with the
$H$-action \cite[Proposition 3.6]{a-Yakimov09p}. The
$H$-prime ideals of $U^-[w]$ are hence in one-to-one correspondence to
the $H$-prime ideals of $R^w_0$ which contain $Q(w)^+_w$. These are
known to be of the form $Q(y,w)_w$ for all $y\in W^{\le w}$ by
\cite{a-Gorelik00} (cf.~\cite[Theorem 3.1]{a-Yakimov09p}). By definition
$P^-(y)=\phi_w(Q(y,w)_w)$ and hence 
\begin{align}\label{eq:YG-equiv}
  U^-[w]/P^-(y)\cong R^w_0/Q(y,w)_w.
\end{align}
Recall that $\omega(U^-[w])=U^+[w]$ where $\omega$ denotes the
involutive algebra isomorphism defined by \eqref{eq:om-def}. For any
$y\le w$ let $P^+(y)=\omega(P^-(y))$. Then $P^+(y)$ is an $H$-prime ideal
of $U^+[w]$ since $\omega $ is an algebra isomorphism and $\omega
(U^0)=U^0$. By the results for $U^-[w]$ explained above,
any $H$-prime ideal of $U^+[w]$ is of the form $P^+(y)$ for some $y\in
W^{\le w}$. Moreover, $\omega$ induces an algebra isomorphism 
\begin{align}\label{eq:-+}
   U^-[w]/P^-(y)\cong  U^+[w]/P^+(y).
\end{align}
Recall the definition of $T^w$ from \eqref{eq:Tw-def}. To identify the
set of polynomial $H$-prime ideals of $U^+[w]$ with a subset of
$W^{\le w}$ define 
\begin{align*}
  W^w=\{w_\Theta\,|\,\Theta\in T^w\}.
\end{align*}
By Lemma \ref{lem:w-sub-red} one has $W^w\subseteq W^{\le w}$. Define
a map  
\begin{align}\label{eq:kappaw}
  \varkappa^w:T^w\rightarrow W^w,\quad  \varkappa^w(\Theta)=w_\Theta.
\end{align}  
The set $T^w$ is partially ordered by inclusion while the set $W^w$ is
partially ordered by the Bruhat order on $W$.
\begin{lem}\label{lem:kappa-bij}
  The map $\varkappa^w$ is an order reversing bijection.
\end{lem}
\begin{proof}
  The map $\varkappa^w$ is surjective by definition and order
  reversing by Lemma \ref{lem:w-sub-red}.
  It follows from \cite[Lemma 2.7.2]{b-BjBr06} that for any $\Theta\in
  T^w$ there exists at most one sequence $j_1<j_2<\dots<j_{|\Theta|}$
  such that the elements $w_m:=s_{j_m}s_{j_{m+1}}\cdots
  s_{j_{|\Theta|}}w$ satisfy
  \begin{align*}
    w_\Theta=w_1<w_2<\dots<w_{|\Theta|}<w.
  \end{align*}
  Hence $\varkappa^w$ is injective by Lemma \ref{lem:w-sub-red} and the fact
  that the reflections $s_\beta$ for $\beta\in \Theta$ commute for any
  $\Theta\in T^w$. 
\end{proof}
\begin{prop}\label{prop:not-poly}
  Let $y\in W^{\le w}\setminus W^w$. Then $U^+[w]/P^+(y)$ is not commutative. 
\end{prop}
\begin{proof}
  By \cite[2.2.6]{b-BjBr06} one can choose $\gamma
  _1,\dots,\gamma _r\in \wurz^+$ such that
  \begin{align}
    s_{\gamma_r}s_{\gamma_{r-1}}\cdots s_{\gamma_1}w=&y,\label{wJ}\\ 
    l(s_{\gamma_i}s_{\gamma_{i-1}}\cdots
    s_{\gamma_1}w)=&l(w)-i\label{w-sequence}
  \end{align}  
   for all $i=1,\dots,r$. As $y\notin W^w$ there exist $a,b\in
   \{1,\dots,r\}$, $a\neq b$, such that $(\gamma_a,\gamma_b)\neq 0$. 
   Now we may apply
   Corollary \ref{cor:12} and assume that $(\gamma_1,\gamma_2)\neq 0$.
   Hence, by Proposition \ref{prop:not-com}, the algebra
   $R^w_0/Q(s_{\gamma_2}s_{\gamma_1}w,w)_w$ is not commutative. By
   \eqref{w-sequence} we have $y\le s_{\gamma_2}s_{\gamma_1}w$ and
   hence $Q(y,w)_w\subseteq Q(s_{\gamma_2}s_{\gamma_1}w,w)_w$. Thus
   $R^w_0/Q(y,w)_w$ is also not commutative which by
   \eqref{eq:YG-equiv} and \eqref{eq:-+} proves the proposition.
\end{proof}
 Recall that for any
$\Theta\in T^w$ we have defined $P_\Theta$ to be the ideal generated by
$\{ E_{\beta}\,|\,\beta\in \wurz ^+_w \setminus \Theta \}$.
Then $P_\Theta\in \mbox{$H-\spec^{poly}(\Uq^+[w])$}$ by
Proposition~\ref{prop:w-reducing}.
\begin{cor}\label{cor:summary}
  (1) The map $W^w\rightarrow \mbox{$H-\spec^{poly}(\Uq^+[w])$}$,
  $y\mapsto P^+(y)$, is an order preserving bijection.\\
  (2) One has $P_\Theta=P^+(w_\Theta)$ for all $\Theta\in T^w$.\\
  (3) If $y\in W^w$ then $U^+[w]/P^+(y)$ is a commutative polynomial
  ring in $l(w)-l(y)$ variables. These variables can be chosen to be the
  $E_{\beta}$ with $\beta\in (\varkappa^w)^{-1}(y)$.
\end{cor}
\begin{proof}
  The map in (1) is order preserving and injective because
  the map $W^{\le w}\rightarrow H-\spec(\Uq^+[w])$, $y\mapsto
  P^+(y)$, is an order preserving bijection. It follows from  
  Proposition \ref{prop:w-reducing} that
  \begin{align}\label{eq:le}
    |T^w|\le | \mbox{$H-\spec^{poly}(\Uq^+[w])$}|.
  \end{align}
  On the other hand Proposition \ref{prop:not-poly}, together with the
  bijection between $W^{\le w}$ and the set of $H$-prime ideals in $U^+[w]$,
  implies that
  \begin{align}\label{eq:ge}
    | \mbox{$H-\spec^{poly}(\Uq^+[w])$}|\le |W^w|.
  \end{align}
  By Lemma \ref{lem:kappa-bij} one obtains $|W^w|=|
  \mbox{$H-\spec^{poly}(\Uq^+[w])$}|$ which shows that the map in (1)
  is bijective. 

  By Proposition \ref{prop:w-reducing} the map $T^w\rightarrow
  H-\spec^{poly}(U^+[w])$, $\Theta\mapsto P_\Theta$ is injective.
  By Lemma \ref{lem:kappa-bij} and (1) we conclude that 
  \begin{align}\label{eq:equal-primes}
    \{P_\Theta\,|\,\Theta\in T^w\}=\{P^+(y)\,|\,y\in W^w\}.
  \end{align}
  To prove (2) it hence suffices to show that $E_{\beta}\notin
  P^+(w_\Theta)$ for all $\Theta\in T^w$ and
  $\beta\in \Theta$. Let $\Theta\in T^w$ and
  $\beta\in \Theta$.  Then $ w_\Theta\le s_{\beta} w$ by Lemma 
  \ref{lem:w-sub-red} and hence $P^+(w_\Theta)\subseteq P^+(s_{\beta}w)$.
  Thus we only need to show that $E_{\beta}\notin P^+(s_{\beta}w)$. To
  this end recall that by definition
  $P^+(s_{\beta}w)=\omega(\phi_w(Q(s_{\beta}w,w)_w))$. Hence
  $\omega(\phi_w(c^{-\lambda}_w c^\lambda_{s_{\beta}w}))\notin
  P^+(s_{\beta} w)$ by Corollary \ref{co:not0Lem}. By definition of
  $\phi_w$ (cf.~\cite[Theorem 3.7]{a-Yakimov09p}) the element
  $\omega(\phi_w(c^{-\lambda}_w c^\lambda_{s_{\beta}w}))$ belongs to
  the weight space of weight $m\beta$ of $U^+[w]$ for some $m\in
  \N_0$. We may assume $m\neq 0$ by choosing $\lambda\in \weights^+$
  regular. Hence $U^+[w]/P^+(s_{\beta} w)$ contains an element of
  weight $m \beta$ for some positive integer $m$.
  
  On the other hand it follows from (1) and
  Equation~\eqref{eq:equal-primes}  that $P^+(w)=P_\emptyset$ and that
  $P^+(s_\beta w)=P_{\{\gamma\}}$ for some $\gamma \in \wurz^+_w$ with
  $\{\gamma\}\in T^w$. Proposition
  \ref{prop:w-reducing} hence implies that $U^+[w]/P^+(s_{\beta} w)$
  is a polynomial ring in $E_{\gamma}$. Therefore $\beta=\gamma$ and
  $E_{\beta}\notin P^+(s_{\beta}w)$ which completes the proof of (2). 
  Property (3) follows immediately from (2) and Proposition
  \ref{prop:w-reducing}. 
\end{proof}
With the above corollary we have all the information necessary to give an
explicit description of the characters of $U^+[w]$. For any two sets $A,B$ let
$\Map(A,B)$ denote the set of maps from $A$ to $B$. As before, let
$\qfield^\times=\field\setminus\{0\}$. Recall also the definition of
$\varkappa^w$ from Equation \eqref{eq:kappaw}.
\begin{thm}\label{thm:char-class}
  There is a bijection
  \begin{align*}
    \Psi:\big\{(y,f)\,|\,y\in W^w, f\in
    \Map((\varkappa^w)^{-1}(y),\qfield^\times)\big\} 
    \rightarrow  \Char(\Uq^+[w])
  \end{align*}
  uniquely determined by
  \begin{align*}
    \Psi(y,f)(E_{\beta})=
    \begin{cases}
      f(\beta)& \mbox{if $\beta \in (\varkappa^w)^{-1}(y)$,}\\
      0 & \mbox{otherwise.} 
    \end{cases}
  \end{align*}
  The inverse map is given by
  $\Psi ^{-1}(\phi )=(w_{\wurz^+_w(\phi)},f_\phi)$ for all $\phi \in
  \Char (\Uq^+[w])$, where 
  $f_\phi(\beta)=\phi (E_{\beta})$ for all $\beta\in \wurz^+_w(\phi )$.
\end{thm}
For the classification of right coideal subalgebras in Theorem
\ref{thm:coid-class} it is necessary to determine subgroups $L$ of the
root lattice which are orthogonal to $\supp\, \phi$ for a given
character $\phi$ of $U^+[w]$. To this end we also note the following
immediate consequence of Proposition \ref{prop:poly}. 
\begin{prop}
  Let $\phi\in \Char(U^+[w])$. Then $\supp\,\phi=\sum_{\beta\in
    \Phi_w^+ (\phi)}\N_0\beta$. 
\end{prop}
\appendix

\section{Root systems and Weyl group combinatorics}
We collect here some auxiliary results about root systems and
the Bruhat order for
finite Weyl groups, which are used to obtain the combinatorial classification of
characters of $U^+[w]$ in Section \ref{characters}.
First we state two lemmata on roots.

\begin{lem} \cite[9.4]{b-Humphreys}
  Let $\al ,\beta \in \wurz $. If $(\al ,\beta )\not=0$ and $(\al ,\al )\le
  (\beta ,\beta )$ then $2|(\al ,\beta )|=(\beta ,\beta )$.
  \label{le:rank2roots}
\end{lem}

\begin{lem} \label{le:nonorth}
Let $\gamma_1,\dots,\gamma_s,\beta,\beta'\in \wurz^+$ and
$a_{1},\dots,a_{s}\in \N _0$. 
Assume that $\beta+\beta'=\sum _{m=1}^{s}a_m\gamma _m$. Then there
exists $n\in \{1,\dots,s\}$ with $a_n(\gamma _n,\beta)\not=0$
or $a_n(\gamma _n,\beta')\not=0$.
\end{lem}

\begin{proof}
This follows from $(\beta +\beta',\beta +\beta')\not=0$.
\end{proof}

Now we turn to
the well-known characterization of the Bruhat order in terms of
positive roots and regular weights. Recall that $\lambda\in \weights^+$ is
called regular if $(\alpha ,\lambda)>0$ for all simple roots $\alpha \in \Pi$.
\begin{lem}\label{le:weyl-bruhat}
  Let $\lambda\in \weights^+$ be a regular weight, $u\in W$, and $\beta\in
  \wurz^+$. Then the following are equivalent:
  (1) $s_\beta u < u$, \qquad  (2) $u^{-1}\beta \in \wurz^-$,\qquad
  (3) $(\beta,u\lambda)<0$.
\end{lem}
\begin{proof}
  The equivalence of (1) and (2) is an immediate consequence of the
  strong exchange condition, cf.~\cite[1.4.3, 4.4.6]{b-BjBr06}.
  The equivalence of (2) and (3) follows from the $W$-invariance of
  the bilinear form $(\cdot,\cdot)$.
\end{proof}
The following technical result is used in the proof of Proposition
\ref{prop:not-poly}. 
\begin{lem}\label{lem:12}
  Let $\alpha,\beta, \gamma\in \wurz^+$ and $w\in W$ such that
  \begin{align}\label{eq:good-length1}
    l(s_\alpha s_\beta s_\gamma w)= l(s_\beta s_\gamma w)-1=l(s_\gamma
    w)-2=l(w)-3.
  \end{align}
  If $(\alpha,\beta)\neq 0$ or $(\alpha,\gamma)\neq 0$ then there
  exist $\alpha',\beta', \gamma'\in \wurz^+$  with
  $(\beta',\gamma')\neq 0$ such that $s_\alpha s_\beta s_\gamma =
  s_{\alpha'}s_{\beta'}s_{\gamma'}$  and
  \begin{align}\label{eq:good-length2}
    l(s_{\alpha'} s_{\beta'} s_{\gamma'} w)= l(s_{\beta'} s_{\gamma'}
    w)-1=l(s_{\gamma'} 
    w)-2=l(w)-3.
  \end{align}
\end{lem}
\begin{proof}
  It suffices to consider the case with $(\beta ,\gamma )=0$.

  If $(\alpha,\beta)=0$ then $s_\alpha s_\beta=s_\beta s_\alpha$ and
  we are done with  $\alpha'=\beta$, $\beta'=\alpha$, and
  $\gamma'=\gamma$. Exchanging $\beta$ and $\gamma$ if necessary we
  may hence assume that 
  \begin{align*}
    (\alpha,\beta)\neq 0 \neq (\alpha,\gamma).
  \end{align*}
  Define $v=s_\gamma w$.\\
  {\bf Case 1: $v^{-1}\alpha<0$.} By Lemma \ref{le:weyl-bruhat} we
  have $l(s_\alpha v)<l(v)$. It remains to show that 
  \begin{align}\label{eq:goal-case1}
    l(s_{s_\alpha\beta}s_\alpha v)<l(s_\alpha v)
  \end{align}
  because then $\alpha'=s_\alpha \beta$, $\beta'=\alpha$, and
  $\gamma'=\gamma$ fulfill the conditions of the lemma. For
  any regular $\lambda\in \weights^+$ one has $(s_\alpha \beta, s_\alpha
  v\lambda)=(\beta,v\lambda)<0$ by Lemma \ref{le:weyl-bruhat} and the
  second equality in \eqref{eq:good-length1}. Again by Lemma
  \ref{le:weyl-bruhat} this proves \eqref{eq:goal-case1}.\\
  {\bf Case 2: $v^{-1}\alpha>0$.} We first claim that
  $(\alpha,\beta)<0$. Indeed, by the first equality in
  \eqref{eq:good-length1} one has for any regular  $\lambda\in
  \weights^+$ the relation 
  \begin{align}\label{hilf}
    0>(\alpha, s_\beta v\lambda)=(\alpha,v\lambda)-\frac{2(\alpha,\beta
    )}{(\beta ,\beta )}(\beta,v\lambda).
  \end{align}
  As $s_\beta v<v$ one has $(\beta,v\lambda)<0$. Using
  $(\alpha,v\lambda)>0$, which holds by assumption, one now gets the 
  desired $(\alpha,\beta)<0$ from \eqref{hilf}.\\
  {\bf Case 2a: $v^{-1}\alpha>0$ and
    $(\alpha,\alpha)=(\beta,\beta)=-2(\alpha,\beta)$.} This is the
  case whenever $\alpha$ and $\beta$ have the same length. As in Case
  1 one verifies that 
  \begin{align}\label{eq:goal2a}
    s_\beta s_{s_\beta \alpha } v < s_{s_\beta \alpha }v< v.
  \end{align}
   Indeed, $(\alpha,s_\beta v\lambda)<0$ for any regular $\lambda\in
   \weights^+$ by \eqref{eq:good-length1} and Lemma
   \ref{le:weyl-bruhat}. Hence $(s_\beta \alpha , v\lambda)<0$ which
   again by Lemma \ref{le:weyl-bruhat} implies $s_{s_\beta \alpha }v<
   v$. To obtain the first relation in \eqref{eq:goal2a} note that
   $s_\beta \alpha =s_\alpha \beta =\alpha+\beta$ and hence 
  \begin{align*}
    s_{s_\beta \alpha }\beta =s_\beta s_\alpha s_\beta \beta =-\alpha.
  \end{align*}
  Now the first relation in \eqref{eq:goal2a} follows from Lemma
  \ref{le:weyl-bruhat} and
  \begin{align*}
     (\beta,s_{s_\beta \alpha } v\lambda)=(-\alpha,v\lambda)<0. 
  \end{align*}
  {\bf Case 2b: $v^{-1}\alpha>0$ and $(\alpha,\alpha)\neq
    (\beta,\beta)$.} Recall that $(\alpha,\beta)<0$. Since $\beta $ and
    $\gamma $ may be exchanged, we can assume that $(\beta,\beta)\neq(\alpha,
  \alpha)\neq(\gamma,\gamma)$ and $(\alpha,\gamma)<0$,
  $(\alpha,\beta)<0$, and $(\beta,\gamma)=0$. This, however, is
  impossible in a finite root system. Indeed, if $(\beta,\beta)=2$ and
  $(\alpha,\alpha)=2d$ with $d=2$ or $d=3$, then $z:=\alpha+\beta+\gamma$ is a
  root. Moreover, this root satisfies $(z,\alpha)=0$ and
  $(z,\beta)=(z,\gamma)\le 0$ by Lemma~\ref{le:rank2roots}.
  This yields the contradiction
  $(z,z)\le 0$. Similarly, one obtains a contradiction if $\beta$ is
  long and $\alpha$ is short by considering $2\alpha+\beta+\gamma$.
\end{proof}
Lemma \ref{lem:w-sub-red} and repeated application of Lemma
\ref{lem:12} imply the following Corollary. 
\begin{cor}\label{cor:12}
  Let $w\in W$ and assume that $\beta_1,\dots,\beta_m\in \wurz^+$ satisfy
  \begin{align*}
    l(s_{\beta_i}s_{\beta_{i-1}}\cdots s_{\beta_1}w)=l(w)-i
  \end{align*}
  for all $i=1,\dots,m$. Assume, moreover, that there exist $i,j\in
  \{1,\dots,m\}$, $i\neq j$, such that $(\beta_i,\beta_j)\neq
  0$. Then there exist $\gamma_1,\dots,\gamma_m\in \wurz^+$ such that
  \begin{align*}
    l(s_{\gamma_i}s_{\gamma_{i-1}}\cdots s_{\gamma_1}w)=&l(w)-i\quad
    \text{for all $i\in \{1,\dots,m\}$,}
    \\
    s_{\gamma_m}s_{\gamma_{m-1}}\cdots s_{\gamma_1}w
    =&s_{\beta_m} s_{\beta_{m-1}} \cdots s_{\beta_1}w, 
  \end{align*} 
  and $(\gamma_1,\gamma_2)\neq 0$.
\end{cor}


\providecommand{\bysame}{\leavevmode\hbox to3em{\hrulefill}\thinspace}
\providecommand{\MR}{\relax\ifhmode\unskip\space\fi MR }
\providecommand{\MRhref}[2]{%
  \href{http://www.ams.org/mathscinet-getitem?mr=#1}{#2}
}
\providecommand{\href}[2]{#2}

\end{document}